\documentclass[11pt]{amsart}

\usepackage{amsmath}
\usepackage{graphicx}
\usepackage{stix}
\usepackage{tikz}
\usepackage{tikz-cd}
\usepackage{hyperref}
\hypersetup{hidelinks}

\newtheorem{mthm}{Theorem} 

\swapnumbers
\newtheorem{theorem}{Theorem}[section]
\newtheorem{lemma}[theorem]{Lemma} 
\newtheorem{corollary}[theorem]{Corollary}

\theoremstyle{definition}
\newtheorem{definition}[theorem]{Definition}
\newtheorem{remark}[theorem]{Remark}

\newtheorem{example}[theorem]{Example}

\makeatletter
\DeclareFontFamily{U}{MnSymbolA}{}
\DeclareFontShape{U}{MnSymbolA}{m}{n}{
    <-6>  MnSymbolA5
   <6-7>  MnSymbolA6
   <7-8>  MnSymbolA7
   <8-9>  MnSymbolA8
   <9-10> MnSymbolA9
  <10-12> MnSymbolA10
  <12->   MnSymbolA12}{}
\DeclareFontShape{U}{MnSymbolA}{b}{n}{
    <-6>  MnSymbolA-Bold5
   <6-7>  MnSymbolA-Bold6
   <7-8>  MnSymbolA-Bold7
   <8-9>  MnSymbolA-Bold8
   <9-10> MnSymbolA-Bold9
  <10-12> MnSymbolA-Bold10
  <12->   MnSymbolA-Bold12}{}
\DeclareSymbolFont{MnSyA}{U}{MnSymbolA}{m}{n}
\SetSymbolFont{MnSyA}{bold}{U}{MnSymbolA}{b}{n}

\DeclareRobustCommand{\overleftharpoon}{\mathpalette{\overarrow@\leftharpoonfill@}}
\DeclareRobustCommand{\overrightharpoon}{\mathpalette{\overarrow@\rightharpoonfill@}}
\def\leftharpoonfill@{\arrowfill@\leftharpoondown\mn@relbar\mn@relbar}
\def\rightharpoonfill@{\arrowfill@\mn@relbar\mn@relbar\rightharpoonup}

\DeclareMathSymbol{\leftharpoondown}{\mathrel}{MnSyA}{'112}
\DeclareMathSymbol{\rightharpoonup}{\mathrel}{MnSyA}{'100}
\DeclareMathSymbol{\mn@relbar}{\mathrel}{MnSyA}{'320}
\makeatother

\newcommand{\Cn}{\mathrm{C}(n)}
\newcommand{{\Csix}}{\mathrm{C}(6)}

\newcommand{\lk}{\mathrm{Link}}
\newcommand{\N}{\Delta}
\newcommand{\diam}{\mathrm{diam}}
\newcommand{\tflat}{\mathbb{E}^2_\mathbf{tri}}
\newcommand{\hflat}{\mathbb{E}^2_\mathbf{hex}}
\newcommand{\wall}{\mathcal E}
\newcommand{\rwall}{\overrightharpoon{W}(\wall)}
\newcommand{\lwall}{\overleftharpoon{W}(\wall)}
\newcommand{\Z}{\mathbf{Z}}

\newcommand{\plength}[1]{|#1|_\mathbf{p}}
\newcommand{\fdist}{d_\mathbf{f}}

\DeclareMathOperator{\Hull}{Hull}

\DeclareMathOperator{\Int}{Int}
\DeclareMathOperator{\Stab}{Stab}

\definecolor{zachcomment}{rgb}{0.55,0.71,0}

\title{Strict $\Csix$ complexes}
\author{Zachary Munro}
\author{Dani Wise}
\begin{document}
    

    


    





    



\begin{abstract}
    We define \emph{strict $\Cn$} small-cancellation complexes, intermediate to $C(n)$ and $C(n+1)$, and we prove groups acting properly cocompactly on a simply-connected strict $\Csix$ complex are hyperbolic relative to a collection of maximal virtually free abelian subgroups of rank $2$. We study geometric walls in a simply-connected strict $\Csix$ complex, and we use them to prove a convex cocompact (cosparse) core theorem for (relatively) quasiconvex subgroups of strict $\Csix$ groups. We provide an example showing the convex cocompact core theorem is false without the strict $\Csix$ assumption. 
\end{abstract}

\maketitle

\section{Introduction}

Dehn solved the word problem for genus $g\geq 2$ orientable surface groups by exploiting the nonpositive curvature of the hyperbolic plane \cite{dehn:unberUnendlicheDiskontinuierlicheGruppen, dehn:transformationDerKurvenAufZweiseitigenFlachen}. The algorithm is referred to as \emph{Dehn's algorithm}. Tartakovskii initiated small-cancellation theory \cite{tartakovskii:solutionWordProblemGroupsReducedBasis}, which was later greatly simplified  by Greendlinger \cite{greendlinger:DehnsAlgorithmForWordProblem, greendlinger:OnDehnsAlgorithmsForConjugacyAndWordProblems}. Lyndon and Schupp further clarified the theory by advancing a diagrammatic viewpoint \cite{lyndon:onDehnsAlgorithm,schupp:OnDehnsAlgorithmAndTheConjugacyProblem}, where small-cancellation conditions impose combinatorial nonpositive curvature. A combinatorial $2$-complex $X$ is $C'(\lambda)$ if for distinct $2$-cell boundaries $\partial R_1$, $\partial R_2$ in the universal cover, we have $\text{Diameter}(\partial R_1\cap \partial R_2)< \lambda\cdot\text{Circumference}(\partial R_1)$, and a group is $C'(\lambda)$ if it is $\pi_1$ of a $C'(\lambda)$ complex. The $C'(\lambda)$ condition is called a \emph{metric} small-cancellation condition because of the inequality above. There are also \emph{nonmetric} conditions $C(n)$ and $T(m)$. A $2$-complex $X$ is $C(n)$ if no $2$-cell $R$ can be surrounded by fewer than $n$ other $2$-cells in the sense that: $\partial R\subset \bigcup_{i=1}^{n-1} \partial R_i$. One imagines prohibiting $R_1, \ldots, R_{n-1}$ cycling entirely around $R$. The $T(m)$ condition prohibits $q$ $2$-cells cycling around a $0$-cell where $2<q<m$. Complexes satisfying $C(n)-T(m)$ for $\frac1n+\frac1m \leq \frac12$ have combinatorial nonpositive curvature properties, which has group-theoretic consequences for $\pi_1$.

Gromov's hyperbolic groups transformed combinatorial group theory, simultaneously generalizing metric small-cancellation groups and fundamental groups of hyperbolic manifolds \cite{gromov:hyperbolicGroups} through a coarse negative curvature condition. Small-cancellation groups are an important source of examples in the theory of $\delta$-hyperbolicity. For example, $C(7)$, $C(5)-T(4)$, and $C(3)-T(7)$ groups are $\delta$-hyperbolic \cite{gerstenShort:smallCancellationTheoryAndAutomaticGroups}, while $C(6)$, $C(4)-T(4)$, and $C(3)-T(6)$  groups are hyperbolic if and only if they do not contain a $2$-dimensional flat \cite{ivanovSchupp:onTheHyperbolicityOfSmallCancellationGroups}. 

Small-cancellation has evolved alongside other directions of research within geometric group theory as well. The theory of $\mathrm{CAT}(0)$ cube complexes has illuminated many problems of separability in several important classes of groups over the past decades. For example, work of Agol and Wise \cite{agol:theVirtualHakenConjecture, wise:cubulatingSmallCancellationGroups} shows that $C'(\frac16)$ groups are subgroup separable. Wise proved that $C'(\frac16)$ groups act properly cocompactly on a $\mathrm{CAT}(0)$ cube complex. The question was raised as to whether nonmetric small-cancellation groups are similarly cubulated \cite{wise:cubulatingSmallCancellationGroups, jankiewicz:cubicalDimensionAndObstructionsToActions}. However, recent work of Munro and Petyt shows there exists a noncocompactly cubulated $C(6)$ group \cite{munroPetyt:coarseObstructionsToCocompactCubulation}. It is still unknown whether all hyperbolic $C(6)$ groups act properly cocompactly on a $\mathrm{CAT}(0)$ cube complex. 

We define the following new small-cancellation condition. A complex $X$ satisfies the \emph{strict $\Cn$} condition if whenever $\partial R\subset \bigcup_{i=1}^n\partial R_i$ with $R_i\neq R$, there is no $1$-cell in both $\partial R\cap \partial R_i$ and $\partial R\cap \partial R_j$ for $i\neq j$. The most important is the strict $\Csix$ condition, which has the following relationship with other nonmetric conditions.
\[
\cdots C(5)\impliedby C(6)\impliedby \text{strict } \Csix \impliedby  C(7) \cdots
\]

Every [finitely presented] group is $\pi_1$ of a [finite] $C(5)$ complex, hence $C(5)$ has no group-theoretical utility. As mentioned above, $C(6)$ yields nonpositive curvature, but not necessarily hyperbolicity since $\Z^2$ is $C(6)$. The following theorem demonstrates how $\Csix$ groups lie between nonpositive and negative curvature.

\begin{mthm}[\ref{thm:relativelyHyperbolic}]
\label{mthm:relativelyHyperbolic}
    Let $G$ act properly cocompactly on a simply-connected strict $\Csix$ complex. Then $G$ is hyperbolic relative to a collection of virtually $\mathbf Z^2$ subgroups. 
\end{mthm}

We also examine wallspace properties of strict $\Csix$ groups. It is intriguing to ask whether strict $\Csix$ groups act properly cocompactly on $\mathrm{CAT}(0)$ cube complexes. We show that strict $\Csix$ groups have features similar to cubulated (relatively) hyperbolic groups. In particular, a simply-connected strict $\Csix$ complex has a collection of geometric walls, similar to the hyperplanes of a cube complex. However, the collection of geometric walls is generally locally infinite and thus fails to turn the $\Csix$ complex into a genuine wallspace. Nonetheless, the geometric walls enable us to show the face-metric on a strict $\Csix$ complex behaves like the combinatorial metric on a $\mathrm{CAT}(0)$ cube complex. The \emph{face-metric $\fdist$} is a metric on the $2$-cells of a complex, where $\fdist(R,R')=n$ if $R=R_1,\dots, R_{n+1}=R'$ is a minimal length sequence of $2$-cells such that $R_i\cap R_{i+1}\neq\emptyset$. Quasiconvex subgroups of cubulated hyperbolic groups have cocompact convex cores \cite{haglund:finiteIndexSubgroupsGraphProducts, sageevWise:coresForQuasiconvexActions}, and we are able to show the analogous theorem holds for hyperbolic strict $\Csix$ groups with the face-metric. 

\begin{mthm}[\ref{cor:convexCore}]
\label{mthm:cocompactCore}
    Let $G$ be a hyperbolic group acting properly cocompactly on a simply-connected strict $\Csix$ complex $X$. Let $H$ be a quasiconvex subgroup of $G$, and let $C$ be a compact subcomplex of $X$. There exists a $\fdist$-convex $H$-cocompact subcomplex $Y$ of $X$ containing $HC$. 
\end{mthm}

A more general theorem holds for all strict $\Csix$ groups. We show that arbitrary strict $\Csix$ groups satisfy a cosparse convex core property for their relatively quasiconvex subgroups. 

\begin{mthm}[\ref{thm:cosparseCore}]
\label{mthm:cosparseCore}
    Suppose $G$ acts properly cocompactly on a simply-connected $\Csix$ complex $X$. Let $H$ be a relatively quasiconvex subgroup of $G$, and let $C$ be a compact subcomplex of $X$. There exists a $\fdist$-convex $H$-cosparse subcomplex $Y$ containing $HC$.
\end{mthm}

We conclude with an example showing the strict $\Csix$ assumption in Theorem~\ref{mthm:cocompactCore} is essential. That is, we exhibit a proper cocompact action of an infinite hyperbolic group $G$ on a simply-connected $2$-complex $X$ and a cocompact quasiconvex subgroup orbit $HC$ such that the $\fdist$ convex hull of $HC$ is not $H$-cocompact. In fact, $H$ is finite and the convex hull of $HC$ is the entire complex $X$. We now briefly discuss the ideas involved in the proof of Theorem~\ref{mthm:relativelyHyperbolic} while summarizing the sections of the paper. 

In Section~\ref{sec:classicalSmallCancellation}, we review the definitions and theorems of classical small-cancellation theory. This includes Greendlinger's Lemma (Lemma~\ref{lem:greendlinger}) for both disc and annular diagrams. Emphasis is placed on definitions which play an important role in the strict small-cancellation theory. For example, hexagonal tilings of the plane, called \emph{honeycombs}, are central objects in the geometry of strict $\Csix$ complexes. We recall the \emph{no missing $i$-shells} property, a form of local convexity introduced by McCammond and Wise \cite{mccammondWise:coherenceLocalyQuasiconvexityAndPerimeterOf2Complexes} to study coherence and local quasiconvexity. We introduce another form of local convexity, the \emph{no missing $i$-complements}, which allows us to prove isometric embeddedness of immersed honeycombs (Lemma~\ref{lem:immersedHoneycombIsometric}). 

In Section~\ref{sec:strictSmallCancellation}, we define the strict $\Cn$ small-cancellation condition and begin to develop its basic properties. We show that coarse intersections of honeycombs in a simply-connected strict $\Csix$ complex are bounded (Theorem~\ref{thm:strictC6isoFlats}). We define \emph{walls} and \emph{wall carriers}, which are tree-like subcomplexes that separate a simply-connected strict $\Csix$ complex into \emph{halfspaces}. This is analogous to hyperplanes in a $\mathrm{CAT}(0)$ cube complex. Since they are not cocompact, these walls fail to turn a strict $\Csix$ complex into a wallspace. We prove some geometric properties, showing walls and their halfspaces are convex (Corollary~\ref{cor:wallFaceConvex}, Lemma~\ref{lem:halfspaceConvex}). 

In Section~\ref{sec:systolization}, we prove Theorem~\ref{mthm:relativelyHyperbolic}. We begin by giving background on \emph{systolic complexes}, which are combinatorial nonpositively curved simplicial complexes. Wise proved that $\Csix$ groups are systolic by showing the dual of a simply-connected $\Csix$ complex is systolic \cite{wise:sixtolicComplexesAndTheirFundamentalGroups}. This was later generalized to graphical $\Csix$ small-cancellation by Osajda--Prytu\l a \cite{osajdaPrytula:classifyingSpacesForFamiliesofSubgroupsForSystolicGroups}. We call this construction \emph{systolization}. We show that systolization induces a bijection between flats and preserves coarse isolation (Theorem~\ref{thm:tricombToHoneycomb}, Corollary~\ref{cor:isoTriIFFisoHex}). Thus, by a systolic version of the Isolated Flats Theorem due to Elsner \cite{elsner:systolicGroupsWithIsolatedFlats}, this implies strict $\Csix$ groups are relatively hyperbolic.

In Section~\ref{sec:convexCores}, we prove Theorems~\ref{mthm:cocompactCore} and \ref{mthm:cosparseCore}. We define \emph{bent wall-segments}, a particular type of $\fdist$-geodesic important in the proofs of Theorems~\ref{mthm:cocompactCore} and \ref{mthm:cosparseCore}. A study of bent wall-segments and halfspaces shows that quasiconvex subcomplexes of a strict $\Csix$ complex are close to their convex hulls (Theorem~\ref{thm:convexHull}). As a corollary, we obtain Theorem~\ref{mthm:cocompactCore}. Next we collect results on relatively hyperbolic groups for the proof of Theorem~\ref{mthm:cosparseCore}. In particular, we make use of isometric embeddedness of relatively quasiconvex subgroups (Lemma~\ref{lem:quasiIsometricEmbedding}) and relative thinness of geodesic triangles (Theorem~\ref{thm:relativelyThinTriangles}). After recalling the definition of \emph{cosparse}, we prove Theorem~\ref{mthm:cosparseCore}. We conclude by providing a counterexample to Theorem~\ref{mthm:cocompactCore} when $X$ is not strict $\Csix$ (Example~\ref{ex:unboundedHull}).

\section{Classical small-cancellation}
\label{sec:classicalSmallCancellation}

\subsection*{\S\ The $\Cn$ condition and diagrams}

    \begin{definition}
        A map between cell complexes is \emph{combinatorial} if $n$-cells are mapped homeomorphically to $n$-cells. Henceforth, all maps are combinatorial. 
        

        A \emph{path} is a map $P\to X$ where $P$ is a $1$-complex homeomorphic to an interval. The path $P\to X$ is a \emph{trivial} if the interval is a point. A \emph{cycle} is a map $C\to X$ where $C$ is a $1$-complex homeomorphic to $S^1$. The \emph{length} of a path or cycle is the number of $1$-cells in $P$ or $C$, respectively.

        Paths $P_1\to X$ and $P_2\to X$ are \emph{equivalent} if there exists an isomorphism $P_1\to P_2$ making the following diagram commute. Similarly, cycles $C_1\to X$ and $C_2\to X$ are \emph{equivalent} if there exists an isomorphism $C_1\to C_2$ making the analogous diagram commute. We do not distinguish between equivalent paths or cycles.

        \begin{center}
            \begin{tikzcd}[row sep=tiny]
                P_1 \arrow[dr]\arrow[dd]& \\
                                        & X\\
                P_2 \arrow[ur]          &
            \end{tikzcd}
        \end{center}

        If $P\to X$ is a nontrivial closed path, then there is an \emph{associated cycle} $C\to X$, where $C$ is the quotient of $P$ identifying its endpoints. Conversely, if $C\to X$ is a cycle, then $P\to X$ is an \emph{associated closed path} if $C\to X$ is associated to $P\to X$. 
    \end{definition}

    \begin{definition}
        A $2$-complex $X$ is \emph{combinatorial} if for each $2$-cell $R$, the attaching map $\partial R\to X^{(1)}$ is a cycle (for some $1$-complex structure on $\partial R$). Henceforth, all complexes are combinatorial $2$-complexes whose attaching maps are immersions.

        If $R$ is a $2$-cell of a complex $X$, then the attaching map $\partial R\to X$ is the \emph{boundary cycle} of $R$. We denote the boundary cycle by $\partial_c R$. A \emph{boundary path} $\partial_p R$ of $R$ is a closed path associated to $\partial_c R$. 
    \end{definition}

    \begin{definition}
        Let $R_1$ and $R_2$ be $2$-cells of $X$. Let $P\to X$ be a path factoring as $P\to \partial R_1\to X$ and $P\to \partial R_2\to X$. The path $P\to X$ is a \emph{piece} if there does not exist an isomorphism $\partial R_1\to \partial R_2$ making the following diagram commute. 

        \begin{center}
		\begin{tikzcd}
			P \arrow[r] \arrow[d] & \partial R_2 \arrow[d] \\
			\partial R_1 \arrow[r] \arrow[ur] & X
		\end{tikzcd}
	\end{center}
    \end{definition}

    \begin{definition} 
        Let $X$ be a complex. Let $d$ denote the \emph{combinatorial metric} on $X^{(0)}$. That is, $d(x,y)$ is the infimal length of a path joining $x$ to $y$.
    
        The \emph{piece-length $\plength{P}$ of a path} $P\to X$ is the infimal $n$ such that $P$ is a concatenation of $n$ pieces. The \emph{piece-length $\plength{C}$ of a cycle} $C\to X$ is the infimal piece-length of a closed path associated to $C$.
    \end{definition}
    
    \begin{definition}[$\Cn$]
        A $2$-complex $X$ is $\Cn$ if $\plength{\partial_c R}\geq n$ for each $2$-cell $R$.
    \end{definition}

    We will use the following example of a simply-connected $\Csix$ complex.

    \begin{definition}
        A \emph{simple honeycomb} is a $2$-complex isomorphic to the tiling of the Euclidean plane by regular hexagons. A \emph{honeycomb} $\hflat$ is a simple honeycomb whose $1$-cells are subdivided (possibly trivially) and referred to as \emph{hexsides}. A \emph{flat} in $X$ is an embedding $\hflat\to X$. 

        Two hexsides in a $2$-cell of $\hflat$ are \emph{equivalent} if they are opposite in a hexagon. This extends to an equivalence relation on hexsides of $\hflat$. A \emph{width $1$ band} is a minimal subcomplex of $\hflat$ containing an equivalence class of hexsides. Two width $1$ bands are \emph{parallel} if they do not share a $2$-cell. A \emph{width $n$ band} is a connected union of $n$ distinct parallel bands. A \emph{flat annulus} is a compact quotient of a band $B$ by a translation of $\hflat$ stabilizing $B$. A flat annulus has \emph{width $n$} if $B$ has width $n$. 
    \end{definition}

    \begin{definition}
        A \emph{disc diagram} $D$ is a compact, contractible $2$-complex with a fixed embedding into $S^2$. The complement of $D$ in $S^2$ is an open $2$-cell $R$. A disc diagram is \emph{trivial} if it is a single point. The \emph{boundary cycle} $\partial_c D$ of a nontrivial disc diagram is the attaching map $\partial R\to D$ of $R$. A \emph{boundary path} $\partial_p D$ is a closed path associated to the boundary cycle. A \emph{disc diagram in $X$} is a map $D\to X$. 
        
        An \emph{annular diagram} $A$ is a compact, connected $2$-complex homotopy equivalent to $S^1$ with a fixed embedding into $S^2$. The complement of $A$ in $S^2$ is two open $2$-cells $R_1$ and $R_2$. An annular diagram has two \emph{boundary cycles} $\partial_{c,1} A$ and $\partial_{c,2} A$, the attaching maps of $R_1$ and $R_2$. An \emph{annular diagram in $X$} is a map $A\to X$. 

        We say $E$ is a \emph{diagram} if $E$ is either a disc or an annular diagram. A \emph{diagram in $X$} is map $E\to X$.
    \end{definition}

\begin{figure}[h]
        \centering
        \includegraphics[width=0.7\textwidth]{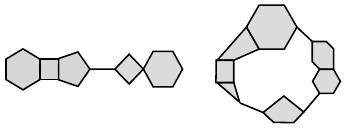}
        \caption{A disc diagram ladder (left) and an annular diagram ladder (right).}	
        \label{fig:ladders}
    \end{figure}

    \begin{definition}
    \label{def:ladder}
        A disc (resp. annular) diagram $E$ is a \emph{ladder} if $E$ is the union of embedded, closed $1$-cells and $2$-cells $R_1,\ldots, R_n$ with $R_i\cap R_j\neq \emptyset$ if and only if $|i-j|\leq 1$ (resp. $|i-j|\leq 1 \mod n$). See Figure~\ref{fig:ladders}.
    \end{definition}

    \begin{definition}
        A diagram $E\to X$ is \emph{reduced} if $P\to E\to X$ is a piece in $X$ for each piece $P\to E$ in $E$.
    \end{definition}

    \begin{definition}
        A \emph{spur} is a degree one $0$-cell in $\partial E$. 
        A $2$-cell $R$ of a diagram $E$ is an \emph{$i$-shell} if $\partial_c R=QS$ where $Q=\partial_c R\cap \partial_cE$ and $\plength S\leq i$. More precisely, $Q$ is a maximal common subpath of $\partial_c R$ and $\partial_c E$. The decomposition $\partial_c R=QS$ is unique, and we call $Q$ and $S$ are the \emph{outerpath} and \emph{innerpath} of $R$. 
    \end{definition}

    Van Kampen's Lemma was discovered in \cite{vankampen:onSomeLemmasInThetheoryOfGroups, lyndon:onDehnsAlgorithm}. See \cite[lem~2.17, lem~2.18]{mccammondWise:fansAndLadders} for the following version.
    
    \begin{lemma}[van Kampen's Lemma]
    \label{lem:vanKampen}
        Let $X$ be a complex. 

        If $C\to X$ is a nullhomotopic cycle, then there exists a reduced disc diagram $\phi:D\to X$ with $\phi\circ(\partial_cD)=C$.
        
        If $C_1\to X$ and $C_2\to X$ are essential, homotopic cycles, then there exists a reduced annular diagram $\phi:A\to X$ with $\phi\circ(\partial_{c,1} A)=C_1$ and $\phi\circ(\partial_{c,2} A)=C_2$.
    \end{lemma}

    By considering the associated cycle, we see that if $P\to X$ is a closed nullhomotopic path, there exists a reduced diagram $\phi:D\to X$ with $\phi\circ(\partial_pD)=P$ for some $\partial_p D$.

    \begin{corollary}
    \label{cor:quotientAnnulus}
        Let $g$ be an automorphism of a simply-connected complex $X$ acting freely and stabilizing subcomplexes $Y_1$ and $Y_2$. There exists a reduced annular diagram $\phi:A\to \langle g\rangle \backslash X$ such that for $i\in \{1,2\}$ there is a map $\partial_iA\to \langle g\rangle \backslash Y_i$ making the following diagram commute.
    \end{corollary}

\begin{figure}[h]
        \centering
        \includegraphics[width=0.4\textwidth]{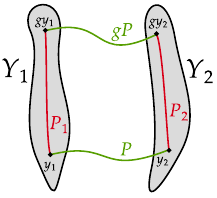}
        \caption{From the proof of Corollary~\ref{cor:quotientAnnulus}.}	
        \label{fig:annulus}
    \end{figure}

    \begin{proof}
        For any points $y_i\in Y_i$ and paths $P_i\to Y_i$ joining $y_i$ to $gy_i$, the cycles $C_i\to \langle g \rangle \backslash Y_i$ associated to $P_i\to \langle g \rangle \backslash Y_i$ are freely homotopic. Indeed, if $P\to X$ is a path joining $y_1$ to $y_2$, then the closed path $PP_2(gP^{-1})P_1^{-1}$ is nullhomotopic. Pushing $P_1\to X$ across such a nullhomotopy to $P_2\to X$ produces a free homotopy of $P_1\to X$ to $P_2\to X$ which descends to a free homotopy of $C_1\to \langle g \rangle \backslash X$ to $C_2\to \langle g \rangle \backslash X$. See Figure~\ref{fig:annulus}. By Lemma~\ref{lem:vanKampen}, there exists a reduced annular diagram $\phi:A\to \langle g \rangle \backslash X$ whose boundary cycles are sent to the $C_i\to \langle g \rangle \backslash Y_i$.
    \end{proof}

    The following is due to Greendlinger. See \cite[thm~V.4.5]{lyndonSchupp:combinatorialGroupTheory}. The version below is a combination of \cite[thm~9.5]{wise:fromRichestoRaags} and \cite[lem~4.2]{bigdelyWise:C(6)GroupsDoNotContainF2xF2}.

    \begin{lemma}[Greendlinger's Lemma]
    \label{lem:greendlinger}
        Let $X$ be a $\Csix$ complex. 
        
        If $D\to X$ is a nontrivial reduced disc diagram in $X$, then either $D$ contains three $3$-shells/spurs, or $D$ is a ladder. In the former case, if $D$ contains exactly three shells/spurs, then the shells must be $2$-shells.

        If $A\to X$ is a reduced annular diagram in $X$, then either $A$ contains a $3$-shell/spur, $A$ is flat, or $A$ is a ladder. 
    \end{lemma}

    \begin{corollary}
    \label{cor:annular2shellOrCycle}
        Let $\phi:A\to X$ be a reduced annular diagram in a $\Csix$ complex, and let $\partial_{c,1} A$ and $\partial_{c,2} A$ be its two boundary cycles. Then either $A$ has a spur, $\partial_{c,1} A=\partial_{c,2} A$, or there exists a $2$-cell $R$ so that $\plength{\phi(\partial_c R\cap \partial_{c,i}A)}\geq 2$ for some $i\in\{1,2\}$.
    \end{corollary}

    \begin{proof}
        A flat annulus of width $1$ is a ladder, and a flat annulus of width $n>1$ always contains a $4$-shell. Thus, by Lemma~\ref{lem:greendlinger}, either $A$ contains a $4$-shell/spur, or $A$ is a ladder. Thus it suffices to consider the cases where $A$ contains a $4$-shell, or $A$ is a ladder.
    
        Suppose $R$ is a $4$-shell, and consider the decomposition $\partial_c R=QS$ where $Q=\partial_c R\cap \partial_{c,i} A$ for some $i\in \{1,2\}$ and $\plength S\leq 4$. Then $\plength Q \geq 2$ otherwise $\plength{\partial_c(\phi R)}\leq \plength{\phi S}+\plength{\phi Q}\leq 5$, contradicting the $\Csix$ condition. 

        If $A$ is a ladder containing no $2$-cells, then $\partial_{c,1} A=\partial_{c,2} A$, so suppose $A$ contains a $2$-cell $R$. The boundary cycle decomposes as $\partial_c R=Q_1S_1Q_2S_2$ where $Q_i=\partial_c R\cap \partial_{c,i} A$ and $\plength{S_i}\leq 1$. Thus $\plength{Q_i}\geq 2$ for some $i\in \{1,2\}$ otherwise $\plength{\partial_c(\phi R)}\leq 4$, a contradiction.
    \end{proof}

\subsection*{\S\ No missing shells and no missing complements}

    The no missing shells property was first introduced in \cite{mccammondWise:coherenceLocalyQuasiconvexityAndPerimeterOf2Complexes}. Here we additionally introduce a related notion, the no missing complements property.

    \begin{definition}
        An immersion $Y\to X$ has \emph{no missing $i$-shells} if the following holds: For any $2$-cell $R\to X$ and decomposition $\partial_c R=QS$ with $\plength S\leq i$, any lift $Q\to Y$ of $Q\to R\to X$ extends to a lift $R\to Y$ of $R\to X$. The immersion $Y\to X$ has \emph{no missing $i$-complements} if for any $2$-cell $R\to X$ and decomposition $\partial_c R=QS$ with $\plength Q\geq i$, any lift $Q\to Y$ extends to a lift $R\to Y$.
    \end{definition}

    Regarding a subcomplex $Y\subset X$ as an inclusion map $Y\to X$, we can also say a subcomplex has no missing $i$-shells. 

    We collect some easy observations about no missing $i$-shells/$i$-complements. If subcomplexes $Y_1,Y_2\subset X$ both have no missing $i$-shells/$i$-complements, then $Y_1\cap Y_2$ has no missing $i$-shells/$i$-complements. For $j>i$, no missing $j$-shells implies no missing $i$-shells, and no missing $i$-complements implies no missing $j$-complements. If $Y\to X$ and $Z\to Y$ have no missing $i$-shells/$i$-complements, then $Z\to Y\to X$ has no missing $i$-shells/$i$-complements. Covering maps have no missing $i$-shells/$i$-complements for all $i$. If $X$ is $\Cn$, then no missing $i$-complements implies no missing $(n-i)$-shells.


    The following lemma is a version of \cite[thm~13.3]{wise:fromRichestoRaags}, which we prove here for completeness. 

    \begin{lemma}
    \label{lem:noMSembedding}
        Let $X$ be a simply-connected $\Csix$ complex. If $Y$ is connected and $Y\to X$ has no missing $3$-shells, then $Y$ is simply-connected and $Y\to X$ is an embedding. Thus we may identify $Y$ with its image $Y\subset X$. 
    \end{lemma}

    \begin{proof}
        If $Y\to X$ has no missing $3$-shells, then so does the composition $\widetilde Y\to Y\to X$. Thus it suffices to show $\widetilde Y\to X$ is an embedding. 
        
        Suppose $\widetilde Y\to X$ is not an embedding. Then there exists an immersed, nonclosed path $P\to \widetilde Y$ with $P\to \widetilde Y\to X$ closed. For each such $P\to \widetilde Y$, there exists a reduced disc diagram $\phi:D\to X$ with $\phi\circ (\partial_p D)=P$ for some boundary path of $D$ by Lemma~\ref{lem:vanKampen}. Choose $P$ and $D$ as above so that $D$ has the minimal number of cells. The diagram $D$ is spurless by minimality. By Lemma~\ref{lem:greendlinger}, there is a $3$-shell $R$ of $D$ whose outerpath $Q$ is a subpath of $P$.

        Let $R\to X$ be the restriction of $\phi$ to $R$. As $Q$ is a subpath of $P\to \widetilde Y\to X$, there is a lift $Q\to \widetilde Y$ of $Q\to R\to X$. Since $\widetilde Y\to X$ has no missing $3$-shells, $Q\to \widetilde Y$ extends to a lift $R\to \widetilde Y$ of $R\to X$. Removing $R$ from $D$ and replacing the outerpath $Q\subset P$ by the innerpath yields a smaller diagram, contradicting minimality.
    \end{proof}

    \begin{lemma}
    \label{lem:noMChelly}
        Let $X$ be a simply-connected $\Csix$ complex, and let $Y_1,\ldots, Y_n$ be connected subcomplexes with no missing $2$-complements. If each $Y_i\cap Y_j\neq \emptyset$, then $\bigcap_iY_i$ is nonempty, connected, and has no missing $2$-complements. 
    \end{lemma}

    \begin{proof} 
        This is trivial in the base case when $n=1$. We will treat the cases $n=2$ and $n=3$ after which we induct.
        
        Let $n=2$. By assumption, $Y_1\cap Y_2$ is nonempty, and $Y_1\cap Y_2$ has no missing $2$-complements since $Y_1$, $Y_2$ have no missing $2$-complements. We now show that $Y_1\cap Y_2$ is connected. Suppose $Y_1\cap Y_2$ has distinct components $C$ and $C'$. By connectivity of $Y_1$, $Y_2$, there are paths $P_1\to Y_1$ and $P_2\to Y_2$ with the same endpoints, joining $C$ to $C'$. By Lemma~\ref{lem:vanKampen}, there is a disc diagram $\phi:D\to X$ with $\phi\circ (\partial_c D)=P_1P_2^{-1}$. Choose $P_1$, $P_2$, and $D$ as above so that $D$ has the minimal number of cells. Observe that $D$ is spurless, since any spur in $D$ can be removed, adjusting $P_1$ and $P_2$ accordingly, contradicting minimality of $D$. Suppose $D$ contains a shell $R$ with $\plength{\partial_c R\cap P_i}\geq 2$ for some $i\in \{1,2\}$. Then $\phi R$ is contained in $Y_i$, since $Y_i$ has no missing $2$-complements. Thus $R$ can be removed from $D$, again adjusting $P_1$ and $P_2$ accordingly and contradicting minimality. By Lemma~\ref{lem:greendlinger}, either $D$ contains three $3$-shells/spurs, or $D$ is a ladder. In the former case, since $D$ contains no spur, it contains three $3$-shells. Thus $D$ contains a $3$-shell $R$ whose outerpath lies in $P_1$ or $P_2$. Such a shell $R$ satisfies $\plength{\partial_c R\cap P_i}\geq 3$ for some $i\in \{1,2\}$, a contradiction. If $D$ is a ladder with $R_1,\ldots, R_n$ as in Definition~\ref{def:ladder}, then $R_1$ is a $1$-shell since $D$ is spurless. The outerpath $Q$ of $R$ satisfies $\plength Q\geq 5$, so $\plength{Q\cap P_i}\geq 2$ for some $i\in\{1,2\}$, a contradiction. Hence $Y_1\cap Y_2$ is connected.

        Let $n=3$. It suffices to show $Y_1\cap Y_2\cap Y_3$ is nonempty, as the remaining properties follow by applying the $n=2$ case to $Y_1\cap Y_2$ and $Y_2\cap Y_3$. Suppose for contradiction that $Y_1\cap Y_2$, $Y_2\cap Y_3$, and $Y_3\cap Y_1$ are pairwise disjoint. Consider a cycle of the form $C=P_1P_2P_3$ where $P_i\to Y_i$ has endpoints in $Y_i\cap Y_{i+1}$ and $Y_i\cap Y_{i+2}$ for each $i\in \Z/3\Z$. Choose $C$ and a diagram $\phi: D\to X$ with $\phi\circ (\partial_c D)=C$ so that $D$ has a minimal number of cells. As above, $D$ is spurless. By Lemma~\ref{lem:greendlinger}, $D$ contains at least two $3$-shells. Thus there exists a shell $R$ whose outerpath $Q$ is contained in $P_iP_{i+1}$ for some $i\in \Z/3\Z$. Since $\plength Q\geq 3$, we have $\plength{Q\cap P_j}\geq 2$ for some $j\in \{i,i+1\}$. As $Y_j$ has no missing $2$-complements, $\phi R$ is a $2$-cell of $Y_j$. Thus, removing $R$ from $D$ and adjusting $P_j$ accordingly yields a smaller diagram, a contradiction.
        
        Suppose the statement holds for $k<n$, and consider $Y_1,\ldots, Y_n$. For $i\in \Z/3\Z$,  $Z_i=\bigcap_{i\neq i+1,i+2}Y_i$ satisfy the hypotheses. Applying the $n=3$ case, $\bigcap_{i\in \Z/3\Z} Z_i=\bigcap_j Y_j$ is nonempty, connected, and has no missing $2$-complements.
    \end{proof}

    \begin{example}
        The no missing $2$-complements condition in Lemma~\ref{lem:noMChelly} is optimal in several ways. 
    
        Consider the $2$-complex $X$ constructed from a $2n$-gon $R$ with sides $P_1,\ldots, P_{2n}$ and bigons $R_1,\ldots, R_{2n}$ such that $R_i$ has one side attached to $P_i$. Evidently, $X$ is a simply-connected $\mathrm{C}(2n)$ complex. Let $Y_1=\bigcup_{i=1}^nR_i$ and $Y_2=\bigcup_{i=n+1}^{2n}R_i$. Both $Y_1$, $Y_2$ are connected subcomplexes of $X$ with no missing $(n-1)$-shells. However, $Y_1\cap Y_2$ is two disconnected points. This shows Lemma~\ref{lem:noMChelly} cannot be strengthened by weakening the assumptions on the $Y_j$ to no missing $i$-shells (for any $i$).    

        A width $1$ band of $\hflat$ has no missing $3$-complements, and it is possible to choose three nonparallel width $1$ bands with empty intersection. This shows that Lemma~\ref{lem:noMChelly} does not hold for no missing $j$-complements with $j>2$.
    \end{example}

    \begin{corollary}
    \label{cor:2cellIntersection}
        If $R_1,\ldots, R_n$ are distinct, pairwise intersecting $2$-cells in a simply-connected $\Csix$ complex, then $\bigcap_i R_i$ is a nonempty, connected subcomplex of $\partial R_1$.
    \end{corollary}

    \begin{proof}
        A closed $2$-cell has no missing $2$-complements. By Lemma~\ref{lem:noMChelly}, $\bigcap_i R_i$ is nonempty and connected. The intersection lies in $\partial R_1$ since $R_1,\ldots, R_n$ are distinct.
    \end{proof}

    \begin{definition}[Faces, face-metric $\fdist$, traces]
    \label{def:faceDistance}
        Let $X$ be a $2$-complex. A $1$-cell is \emph{isolated} if it does not lie on a $2$-cell. A \emph{face} of $X$ is a $2$-cell or isolated $1$-cell of $X$. We give the set of faces a metric. The \emph{face-metric} $\fdist$ is defined by letting $\fdist(R, R')$ be the infimal $n$ such that there exists a sequence of faces $R=R_1,\ldots, R_{n+1}=R'$ with $R_i\cap R_{i+1}\neq \emptyset$ for all $i$. We call such a sequence $R_1,\ldots, R_{n+1}$ a \emph{$\fdist$-geodesic}. A path $P\to X$ \emph{traces} a $\fdist$-geodesic $R_1,\ldots, R_m$ if there is a decomposition $P=P_1\cdots P_m$ with $P_i\to R_i$.
    
        Let $Y\subset X$ be a subcomplex. We say $Y$ is \emph{trace-convex} if there exists a trace $P\to Y$ for every $\fdist$-geodesic joining faces of $Y$. We say $Y$ is \emph{face-convex} if every $\fdist^X$-geodesic joining faces of $Y$ is contained in $Y$. 
    \end{definition}


    \begin{lemma}
    \label{lem:noMStraceconvex}
        Let $X$ be a simply-connected $\Csix$ complex. If a subcomplex $Y\subset X$ has no missing $3$-shells, then $Y$ is trace-convex.
    \end{lemma}

    \begin{proof}
        Let $P=P_1\cdots P_{n}$ be a trace of a $\fdist^Y$-geodesic $R_1,\ldots, R_n$, and let $J=J_1\cdots J_m$ be the trace of a $\fdist^X$-geodesic $C_1,\ldots, C_m$ with the same endpoints $C_1=R_1$ and $C_m=R_n$ as $P$. By Lemma~\ref{lem:greendlinger} there exists a disc diagram $\phi:D\to X$ with $\phi\circ (\partial_c D)=PJ^{-1}$. Choose $P$, $J$, and $D$ so that $D$ has the minimal number of cells. 

        By minimality, $D$ does not have a spur in the interior of either $P$ or $J$. 
        
        Suppose $R$ is a $3$-shell of $D$ with outerpath $Q$ and innerpath $S$. The outerpath cannot lie entirely in $P$ or $J$. Indeed, if $Q$ lies in $J$, then $Q$ must intersect at least three consecutive cells $C_{i-1}, C_i, C_{i+1}$. Thus the sequence $C_1,\ldots, C_{i_1}, C, C_{i+1},\ldots, C_m$ is a $\fdist^X$-geodesic, so $Q$ can be replaced by $S$ in $J$ and $R$ can be removed from $D$, contradicting minimality. If $Q$ is a subpath of $P$, then $Q$ intersects a subpath of the form $P_{i-1}, P_i, P_{i+1}$. Additionally, $R$ is a $2$-cell of $Y$, since $Y$ has no missing $3$-shells. Thus $Q$ can be replaced by $S$ in $P$, and $R$ can be removed from $D$, a contradiction. We have deduced that $D$ can have at most two $3$-shells, so by Lemma~\ref{lem:greendlinger}, $D$ is a ladder $L_1,\ldots, L_k$. 

        Let $R$ be the first $L_i$ which is a $2$-cell. Then $\partial_cL=Q_1Q_2Q_3$ where $Q_1=\partial_cL\cap J$, $Q_2=\partial_cL\cap L_{i+1}$, and $Q_3=\partial_cL\cap P$. A minimality argument as above shows that $Q_1$ intersects at most two consecutive $R_i, R_{i+1}$. Thus $\plength{Q_1Q_2}\leq 3$, and $Q_3$ is the outerpath of a $3$-shell contained in $Y$. Sincee $Y$ has no missing $3$-shells, an argument as above leads to a contradiction of minimality. Hence $D$ is a ladder of $1$-cells and $P=J$.
    \end{proof}

    \begin{lemma}
    \label{lem:immersedHoneycombIsometric}
        Let $X$ be a simply-connected $\Csix$ complex. Any immersion $\hflat~\to~X$ is $\fdist$-isometric flat with no missing $3$-complements.
    \end{lemma}

    \begin{proof}
        Let $R$ be a $2$-cell with $\partial_cR=QS$ and $\plength Q\geq 3$. Suppse there exists a lift $Q\to \hflat$. The piece-length of $Q$ in $\hflat$ upper bounds the piece-length of $Q$ in $X$. Thus $Q\to \hflat$ contains two hexsides of a hexagon $R'$ in $\hflat$. Hence $R= R'$, otherwise $\plength{\partial_c R'}\leq 5$, contradicting the $\Csix$ condition.

        We have shown $\hflat\to X$ has no missing $3$-complements (hence no missing $3$-shells), so $\hflat \hookrightarrow X$ by Lemma~\ref{lem:noMSembedding}. For two $2$-cells $R$, $R'$ of $\hflat$, there is a trace $P=P_1\cdots P_n\to \hflat$ of a $\fdist^X$-geodesic joining $R$ and $R'$ by Lemma~\ref{lem:noMStraceconvex}. No $P_i$ contains two hexsides, as this would violate $\Csix$. Hence each $P_i$ lies in a hexagon of $\hflat$. Thus there exists a $\fdist^X$-geodesic joining $R$, $R'$ consisting of hexagons in~$\hflat$.  
    \end{proof}
    
\section{Strict small-cancellation}
\label{sec:strictSmallCancellation}

\subsection*{\S Strict $\Csix$ and isolated flats}

    \begin{definition}[Strict $\Cn$]
        A $2$-complex $X$ is \emph{strict $\Cn$} if $\plength P>n$ for any $2$-cell $R$ and immersed path $P\to \partial R$ with $|P|>|\partial_c R|$. The $\Cn$ condition is equivalent to: $\plength P \geq n$ for any $2$-cell $R$ and immersed $P\to \partial R$ with $|P|\geq |\partial_c R|$. Thus strict $\Cn$ implies $\Cn$.
    \end{definition}

    \begin{definition}[Petals]
        Let $R$ be a $2$-cell of a $\Csix$ complex $X$. Suppose $R_1,\ldots , R_6$ are $2$-cells distinct from $R$ and $\partial_c R=P_1\cdots P_6$ is a decomposition of $\partial_c R$ into pieces with $P_i$ contained in $R\cap R_i$. Each $R_i$ is a \emph{petal} of $R$ with \emph{petal-piece} $P_i$ in $R$.
    \end{definition}

    Each hexside of a $2$-cell $R$ of a honeycomb is a petal-piece in $R$.

    \begin{lemma}
    \label{lem:petalPieceMaximal}
        A $\Csix$ complex is strict $\Csix$ if and only if for any $2$-cell $R$, petal-pieces are maximal pieces in $R$.
    \end{lemma}

    In other words, if $P\to R$ properly contains a petal-piece of $R$, then $\plength P\geq 2$.

    \begin{proof}
        Suppose $R$ contains a nonmaximal petal-piece. Then there exists a petal-piece decomposition $\partial_c R=P_1\cdots P_6$ and a piece $P'_6\to R$ properly containing $P_6$ as an initial subpath. The path $P_1\cdots P_5P'_6\to R$ is longer than $\partial_c R$ and has piece-length six, contradicting strict $\Csix$.

        Conversely, suppose a path $P\to \partial R$ violates strict $\Csix$. Then there exists a piece decomposition $P=P_1\cdots P_5P'_6$, where $P'_6$ can be shortened to $P_6$ so that $P_1\cdots P_6=\partial_c R$ is a petal-piece decomposition.
    \end{proof}

    \begin{lemma}
    \label{lem:noMoreThanHexside}
        Let $\hflat\to X$ be a flat in a simply-connected, strict $\Csix$ complex. If a $2$-cell $R\cap \hflat$ properly contains a hexside of $\hflat$, then $R\subset \hflat$.
    \end{lemma}

    \begin{proof}
        Suppose $R\cap \hflat$ properly contains a hexside of $R'$. Then there is a path $P\to R\cap R'$ properly containing a hexside. If $R\neq R'$, then $P$ is a piece. This contradicts Lemma~\ref{lem:petalPieceMaximal}.
    \end{proof}

    \begin{corollary}
    \label{cor:noMissingPetalPieces}
        A honeycomb $\hflat$ in a strict $\Csix$ complex has no missing $2$-complements. That is, if $\plength{\partial_c R\cap \hflat}\geq 2$ for some $2$-cell $R$, then $R\subset \hflat$.  
    \end{corollary}

    \begin{proof}
        By Lemma~\ref{lem:petalPieceMaximal}, $\plength{\partial_c R\cap \hflat}\geq 2$ if and only if $R$ properly contains a hexside of $\hflat$. By Lemma~\ref{lem:noMoreThanHexside}, $R\subset \hflat$.
    \end{proof}

    \begin{corollary}
    \label{cor:honeycombIntersection}
        Let $\hflat$ and $\mathbb{E'}^2_\mathbf{hex}$ be distinct honeycombs in a simply-connected, strict $\Csix$ complex $X$. Then either $\hflat\cap \mathbb{E'}^2_\mathbf{hex}$ is a single $2$-cell, or $\hflat\cap \mathbb{E'}^2_\mathbf{hex}$ does not properly contain a hexside of either honeycomb.
    \end{corollary}

\begin{figure}[h]
        \centering
        \includegraphics[width=0.4\textwidth]{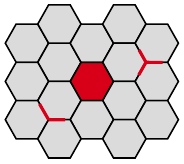}
        \caption{Each red component indicates a possible intersection with another honeycomb.}	
        \label{fig:intersection}
    \end{figure}

    $\hflat \cap \mathbb{E'}^2_\mathbf{hex}$ is connected by Corollary~\ref{cor:noMissingPetalPieces} and Lemma~\ref{lem:noMChelly}. Thus, Corollary~\ref{cor:honeycombIntersection} implies $\hflat \cap \mathbb{E'}^2_\mathbf{hex}$ is either a $2$-cell or a (possibly degenerate) tripod with legs properly contained in hexsides of both $\hflat$ and $\mathbb{E'}^2_\mathbf{hex}$. See Figure~\ref{fig:intersection}.

    \begin{proof}
        If $\hflat\cap \mathbb{E'}^2_\mathbf{hex}$ does not contain $2$-cells, then it does not properly contain a hexside of either $\hflat$ or $\mathbb{E'}^2_\mathbf{hex}$ by Corollary~\ref{cor:noMissingPetalPieces}. Suppose $\hflat\cap \mathbb{E'}^2_\mathbf{hex}$ properly contains a $2$-cell $R_1$. By Corollary~\ref{cor:noMissingPetalPieces} and Lemma~\ref{lem:noMChelly}, $\hflat\cap \mathbb{E'}^2_\mathbf{hex}$  is connected and has no missing $2$-complements. By connectivity, the intersection contains a hexside of a $2$-cell $R_2\subset \hflat$ sharing a hexside with $R_1$. Since $\hflat\cap \mathbb{E'}^2_\mathbf{hex}$ has no missing $2$-complements, it contains $R_2$. There is an enumeration of the $2$-cells $R_1, R_2,\ldots$ of $\hflat$ so that $R_1\cup\cdots\cup R_i$ properly contains a hexside of $R_{i+1}$. By induction, $\hflat\cap \mathbb{E'}^2_\mathbf{hex}$ contains $R_i$ for every $i$, so $\hflat= \hflat\cap \mathbb{E'}^2_\mathbf{hex}$. Likewise, $\mathbb{E'}^2_\mathbf{hex}= \hflat\cap \mathbb{E'}^2_\mathbf{hex}$.
    \end{proof}

    Corollary~\ref{cor:honeycombIntersection} shows that honeycombs are equal when they share two $2$-cells.

    \begin{lemma}
    \label{lem:honeycombsLocallyFinite}
        Let $X$ be a $\fdist$-locally-finite, strict $\Csix$ complex. Then the set of honeycombs in $X$ is locally finite, i.e., for any compact $Y\subset X$, only finitely-many honeycombs intersect $Y$. 
    \end{lemma}

    \begin{proof}
        Let $Y'=N_2(Y)$ be the $2$-neighborhood of $Y$ with respect to the face-metric. By assumption, $Y'$ is finite. A honeycomb intersecting $Y$ must intersect $Y'$ in at least two $2$-cells. By Corollary~\ref{cor:honeycombIntersection}, distinct honeycombs intersecting $Y$ must have distinct intersections with $Y'$, so the lemma follows. 
    \end{proof}

    \begin{theorem}
    \label{thm:honeycombStabilizerCocompact}
        Suppose $G$ acts cocompactly on a simply-connected, strict $\Csix$ complex $X$. Then stabilizers of honeycombs act cocompactly. 
    \end{theorem}

    \begin{proof}
        Let $\hflat$ be a honeycomb in $X$. Consider pairs of adjacent $2$-cells in a band $B$ of $\hflat$. By cocompactness, there are distinct ordered pairs $(R_1$, $R'_1)$ and $(R_2$, $R'_2)$ in $B$ with $gR_1=R_2$, $gR'_1=R'_2$ for some $g\in G$. The flats $\hflat$ and $g\hflat$ share two $2$-cells. By Corollary~\ref{cor:honeycombIntersection}, $\hflat=g\hflat$, and $g$ is a translation stabilizing $B$. Similarly, there is a $g'$ translating $\hflat$ along a band $B'$ nonparallel to $B$. Hence $\hflat$ is $\langle g, g'\rangle$-cocompact.  
    \end{proof}

    \begin{definition}
        A $\Csix$ complex $X$ has \emph{isolated flats} if for every $r>0$ there exists $K>0$ so that $\diam\big(N_r(\hflat)\cap N_r(\mathbb{E'}^2_\mathbf{hex})\big)<K$ for honeycombs $\hflat\neq \mathbb{E'}^2_\mathbf{hex}$, where the diameter is measured with respect to the face-metric. 
    \end{definition}

    If $X$ is a simply-connected cocompact $\Csix$ complex, then the above diameter can equivalently be taken with respect to the combinatorial metric. 

    \begin{theorem}
    \label{thm:strictC6isoFlats}
        Suppose $G$ acts properly cocompactly on a simply-connected, strict $\Csix$ complex $X$. Then $X$ has isolated flats. 
    \end{theorem}

    \begin{proof}
        Suppose $\diam\big(N_r(\hflat)\cap N_r(\mathbb{E'}^2_\mathbf{hex})\big)$ is unbounded as honeycombs vary. By cocompactness, there exists a compact subspace $Y$ so that $\diam\big(N_r(\hflat)\cap N_r(\mathbb{E'}^2_\mathbf{hex})\big)$ is unbounded across pairs of honeycombs intersecting $Y$. By Lemma~\ref{lem:honeycombsLocallyFinite}, there exists a single pair of honeycombs with $\diam\big(N_r(\hflat)\cap N_r(\mathbb{E'}^2_\mathbf{hex})\big)$ unbounded. Let $Y_1=\hflat$ and $Y_2=\mathbb{E'}^2_\mathbf{hex}$ be such a pair. 

        Applying Theorem~\ref{thm:honeycombStabilizerCocompact}, there exists a compact $Y$ and infinitely-many $g_i\in G$ such that $g_i Y_1$, $g_i Y_2$ both intersect $Y$. Applying Lemma~\ref{lem:honeycombsLocallyFinite} again, there exists $g=g_j^{-1}g_i\neq 1$ acting by a nontrivial translation on both $Y_1$ and $Y_2$. By Corollary~\ref{cor:quotientAnnulus}, there exists a reduced annular diagram $\phi: A\to \langle g \rangle \backslash X$ such that $\phi\circ (\partial_{c,1} A)$ and $\phi\circ (\partial_{c,2} A)$ are contained in $Y_1/g$ and $Y_2/g$. Let $A$ be such a diagram with the minimal number of cells. 

        By minimality, $A$ must be spurless. If $R$ is a $2$-cell of $A$ such that $\plength{\partial_c R\cap \partial_{c,i} A}\geq 2$, then a lift of $R$ to $X$ intersects $Y_i$ along a path of piece-length two. By Lemma~\ref{cor:noMissingPetalPieces}, $R$ is contained in $\hflat$, contradicting minimality of $A$. By Corollary~\ref{cor:annular2shellOrCycle}, $A$ is a cycle. Thus $\hflat\cap \mathbb{E'}^2_\mathbf{hex}$ is unbounded, contradicting Corollary~\ref{cor:honeycombIntersection}.
    \end{proof}

\subsection*{\S Walls in strict $\Csix$ complexes}

    \begin{definition}
        Let $X$ be a $2$-complex, and let $R$ be a $2$-cell of $X$. Two $1$-cells $e$ and $e'$ of $R$ are \emph{opposite} if $\plength P >3$ for any path $P\to \partial R$ traversing $e$ and $e'$.
    \end{definition}
    
    \begin{lemma}
    \label{lem:unique2Cell}
        Let $X$ be a simply-connected $\Csix$ complex. Let $e$ and $e'$ be opposite $1$-cells of a $2$-cell $R$. Then $\partial_c R=\partial_c R'$ for any $2$-cell $R'$ containing $e$ and $e'$.
    \end{lemma}

    Note that distinct $2$-cells can have the same boundary cycle. For example, consider the universal covers of the standard $2$-complexes associated with presentations of the form $\langle a\mid a^n\rangle$.

    \begin{proof}
        Suppose $2$-cells $R$ and $R'$ have distinct boundary cycles which both contain $e$ and $e'$. By Corollary~\ref{cor:2cellIntersection}, there exists a path $P\to R\cap R'$ traversing $e$ and $e'$. But then $\plength{P}=1$, a contradiction.
    \end{proof}

    As a consequence of Lemma~\ref{lem:unique2Cell}, we speak of $1$-cells $e$, $e'$ as being opposite without reference to a $2$-cell, since $e$, $e'$ are opposite in any $2$-cell containing $e$, $e'$. 

    \begin{lemma}
    \label{lem:oppositesExist}
        A $2$-complex $X$ is strict $\Csix$ if and only if for every $2$-cell $R$ and $1$-cell $e$ of $R$, there exists a $1$-cell $e'$ of $R$ opposite $e$.
    \end{lemma}

    \begin{proof}
        If $X$ is not strict $\Csix$, then there exists a $2$-cell $R$ and an immersed path $P\to \partial R$ satisfying $\plength{P}\leq 6$ and $|P|>|\partial_c R|$. Let $P=P_1P_2$, where $\plength{P_1}=\plength{P_2}= 3$. Let $e$ be a $1$-cell that $P$ traverses twice. Then $P_1$ and $P_2^{-1}$ have piece-length three, traverse $e$ in opposite directions, and cover $\partial R$. Thus $e$ has no opposite in $R$.  
        
        Suppose $X$ is strict $\Csix$. Let $R$ be a $2$-cell of $X$, and let $e$ be a $1$-cell of $R$. Let $P\to \partial R$ and $J\to \partial R$ be maximal-length paths that begin at $e$, traverse $\partial R$ in opposite directions, and satisfy $\plength{P},\plength{J}\leq 3$. If $\partial R$ is covered by $P$ and $J$, then after removing (possibly trivial) terminal subpaths of $P$ and $J$, the concatenation $PJ^{-1}\to \partial R$ is an immersed path satisfying $\plength{PJ^{-1}}\leq 6$ and $|PJ^{-1}|>|\partial_c R|$. Thus there exists an $1$-cell $e'$ of $R$ not contained in $P$ or $J$, and $e'$ is opposite to~$e$.  
    \end{proof}

    \begin{definition}
        Let $X$ be a simply-connected strict $\Csix$ complex. A collection of $1$-cells $\wall$ is a \emph{semi-wall} if it satisfies the following two conditions:
        \begin{itemize}
            \item For all $e,e'\in \wall$ there is a sequence $e=e_1,\ldots, e_n=e'$ in $\wall$ where $e_i$ and $e_{i+1}$ are opposite for each $i$.
            \item For each $2$-cell $R$ of $X$, either $\wall$ contains no  $1$-cell, one $1$-cell, or exactly two opposite $1$-cells of $R$.
        \end{itemize}
        
        If for each $R$, either $\wall$ contains no $1$-cell or exactly two opposite $1$-cells of $R$, then $\wall$ is a \emph{wall}. The \emph{carrier} $W(\wall)$ of $\wall$ is the union of all faces containing a $1$-cell of $\wall$. The \emph{interior} $\Int(W(\wall))$ is the union of open cells whose closure contains a $1$-cell of $\wall$.
    \end{definition}

    \begin{definition}
        Let $X$ be a simply-connected $\Csix$ complex. Let $R_1,\ldots, R_n$ be a sequence of $2$-cells, and let $e_1,\ldots, e_{n-1}$ be a sequence of $1$-cells such that $e_i$ and $e_{i+1}$ are opposite in $R_{i+1}$. If $\partial_cR_i\neq \partial_cR_{i+1}$ for each $i$, then $R_1,\ldots, R_n$ is a \emph{wall-segment}. 
    \end{definition}

    \begin{lemma}
    \label{lem:wallSegment}
        Let $X$ be a simply-connected $\Csix$ complex. Then a wall-segment $R_1,\ldots, R_n$ is the unique $\fdist$-geodesic joining $R_1$ to $R_n$.
    \end{lemma}

    \begin{proof}
        Suppose the statement is false, and let $R_1,\ldots, R_n$ be a minimal counterexample. That is, there exists a $\fdist$-geodesic $S_1,\ldots, S_m$ with $S_1=R_1$, $S_m=R_n$, and $S_i\neq R_j$ for all other $i,j$. Let $Y$ be the subcomplex $R_1\cup\cdots \cup R_n$. Consider a path $J=J_2\cdots J_{n-1}\to Y$ with $J_i\to R_i$ and a trace $P=P_1\cdots P_m$ of $S_1,\ldots, S_m$ such that $P$ and $J$ have the same endpoints. Choose $P$, $J$, and a disc diagram $D\to X$ with $\partial_c D=PJ^{-1}$ such that $D$ has the minimal number of cells.
        
        Observe that $\plength{J_i}>1$, since $J_i$ joins a piece in $R_{i-1}\cap R_i$ containing $e_{i-1}$ to a piece in $R_i\cap R_{i+1}$ containing $e_i$. In particular, if a $2$-cell $R$ contains $J_i$, then $R=R_i$. Thus, by minimality, no $2$-cell of $D$ can overlap $J$ along a path of piece-length greater or equal to three. Likewise for $2$-cell overlaps with $P$. Also by minimality, $D$ has no spurs in the interiors of either $P$ or $J$. Thus by Lemma~\ref{lem:greendlinger}, $D$ is a ladder $L_1, \ldots, L_k$. We claim $D$ is a ladder of $1$-cells, and so $P=J$. Indeed, let $L_i$ be the first $2$-cell in the ladder. From our above observations, we have $\plength{\partial_c L_i\cap P}<3$ and $\plength{\partial_c L_i\cap J}<3$, but this is impossible by $\Csix$.  

        Lastly, no $P_k$ contains a $J_i$, since $\plength{J_i}>1$. Hence $m>n$, contradicting that $S_1,\ldots, S_m$ is a $\fdist$-geodesic.
    \end{proof}

    \begin{lemma}
    \label{lem:wallSegmentsInWalls}
        Let $X$ be a simply-connected $\Csix$ complex. Any $2$-cells $R$, $R'$ of a semi-wall carrier $W(\wall)$ are joined by a wall-segment $R=R_1,\ldots, R_n=R'$ in~$W(\wall)$.
    \end{lemma}

    \begin{proof}
        Let $e_1,\ldots, e_{n-1}$ be a minimal length sequence of $1$-cells in $\wall$ where $e_1$ is a $1$-cell of $R$, and $e_{n-1}$ is a $1$-cell of $R'$, and $e_i$, $e_{i+1}$ are opposite. Let $R_1=R$, $R_n=R'$, and for each other $i$, let $R_i$ be a $2$-cell containing $e_{i-1}$, $e_i$. By the definition of a semi-wall, if $\partial_c R_i= \partial_c R_{i+1}$ then $e_i=e_{i+2}$, contradicting minimality. Thus $\partial_c R_i\neq \partial_c R_{i+1}$ for each $i$, and the sequence $R=R_1,\ldots, R_n=R$ is a wall-segment by Lemma~\ref{lem:wallSegment}.
    \end{proof}

    \begin{corollary}
    \label{cor:wallFaceConvex}
        Semi-wall carriers in simply-connected $\Csix$ complexes are face-convex. 
    \end{corollary}

    \begin{proof}
        Let $R$ and $R'$ be two faces of a carrier $W(\wall)$. By Lemma~\ref{lem:wallSegmentsInWalls}, there exists a wall segment $R_1,\ldots, R_n$ in $W(\wall)$ where $R_1=R$ and $R_n=R'$. By Lemma~\ref{lem:wallSegment}, $R_1,\ldots, R_n$ is the unique $\fdist$-geodesic joining $R$ and $R'$. In particular, every $\fdist$-geodesic joining $R$ and $R'$ is contained in $W(\wall)$.
    \end{proof}

    \begin{corollary}
    \label{cor:intersecting2CellsInWall}
        Let $X$ be a simply-connected strict $\Csix$ complex. If $2$-cells $R$ and $R'$ of a semi-wall carrier $W(\wall)$ intersect, then $R\cap R'$ contains a $1$-cell of $\wall$. 
    \end{corollary}

    \begin{proof}
        By Lemma~\ref{lem:wallSegmentsInWalls} there exists a wall-segment $R=R_1,\ldots, R_n=R'$ in $W(\wall)$. Since wall-segments are $\fdist$-geodesics and $R\cap R'\neq \emptyset$, we must have $n=2$. Then $R\cap R'$ contains a $1$-cell of $\wall$.
    \end{proof}

    \begin{lemma}
    \label{lem:extendingWall}
        Let $X$ be a simply-connected strict $\Csix$ complex. Any pair of opposite $1$-cells $e$, $e'$ is contained in a wall.
    \end{lemma}

    \begin{proof}
        By Zorn's lemma, there exists a maximal semi-wall $\wall$ containing $e$, $e'$.

        We claim $\wall$ is a wall. Indeed, suppose $R\cap \wall$ is a single $1$-cell $f$ of $R$. By Lemma~\ref{lem:oppositesExist}, there exists a $1$-cell $f'$ of $R$ opposite to $f$. By maximality, $f'$ cannot be added to $\wall$. Thus there exists a $2$-cell $R'$ containing $f'$ such that $R'\cap \wall$ contains two $1$-cells. Corollary~\ref{cor:intersecting2CellsInWall} implies $R\cap R'$ contains a $1$-cell of $\wall$, which must be distinct from $f$ since $R\cap R'$ contains the $1$-cell $f'$ opposite $f$. This contradicts that $R$ contains only one $1$-cell of $\wall$.
        
    \end{proof}

    \begin{definition}
        In light of Lemma~\ref{lem:extendingWall}, given opposite $1$-cells $e$ and $e'$, we let $W(e,e')$ denote a wall carrier containing $e$ and $e'$. 
    \end{definition}

    \begin{lemma}
    \label{lem:wallSeparates}
        Let $X$ be a simply-connected $\Csix$ complex. For any wall $\wall$, the complement $X-\Int(W(\wall))$ has exactly two components.
    \end{lemma}

    \begin{proof}[Sketch]
        There are two opposite $1$-cells $e_1$, $e_2$ in $R\cap \wall$ for each $2$-cell $R$ of $W(\wall)$. For each such triple $e_1$, $e_2$, $R$, let $m_R$ be a closed segment embedded in $R$ whose endpoints are the midpoints of $e_1$ and $e_2$. Let $T=\bigcup_{R\in W(\wall)}m_R$. By Lemma~\ref{lem:wallSegmentsInWalls}, the graph $T$ has no embedded cycles besides bigons $m_Rm_{R'}^{-1}$ where $\partial_cR=\partial_cR'$. Thus $T$ is a tree with multi-edges and has a neighborhood homeomorphic to $T\times [0,1]$. A Mayer-Vietoris sequence shows that $X-T$ has two components. See \cite[lem~2.10]{hsuWise:cubulatingGraphsOfFreeGroupsWithCyclicEdgeGroups}, for example. Consequently, $W(\wall)-T$ also has two components. The set $\Int(W(\wall))$ is an open neighborhood of $T$ which does not disconnect either of the components in $W(\wall)-T$. Thus $X-\Int(W(\wall))$ has two components.
    \end{proof}

    \begin{definition}
    \label{def:halfspace}
        Let $\wall$ be a wall of $X$. By Lemma~\ref{lem:wallSeparates}, there is a decomposition $A\sqcup B=X-\Int(W(\wall))$ into connected components. The subspaces $\rwall = A\cup W(\wall)$ and $\lwall = B\cup W(\wall)$ are the \emph{halfspaces} of $\wall$.
    \end{definition}

    \begin{lemma}
    \label{lem:halfspaceConvex}
        Let $X$ be a simply-connected strict $\Csix$ complex. The halfspaces $\lwall$, $\rwall$ are face-convex for any wall $\wall$ of $X$.
    \end{lemma}

    \begin{proof}
        We wish to show a $\fdist$-geodesic $R_1,\ldots, R_n$ is contained in $\lwall$ if $R_1,R_n \subset \lwall$. Let $\lwall=A\cup W(\wall)$ be the union from Definition~\ref{def:halfspace}. By Lemma~\ref{lem:wallSeparates}, $R_1,\ldots, R_n$ is a concatenation of subsequences contained in $A$ and subsequences starting and ending in $W(\wall)$. By Corollary~\ref{cor:wallFaceConvex}, these latter subsequences lie in $W(\wall)$. Hence $R_1,\ldots, R_n$ is contained in $\lwall$.
    \end{proof}

\section{Systolization}
\label{sec:systolization}

In the context of group theory, systolic complexes were first studied independently by Haglund and Jankiewicz--Swi\c antkowski \cite{januszkiewiczSwiatkowski:simplicialNonpositiveCurvature, haglund:complexesSimpliciauxHyperboliquesDeGrandeDimension}. The $1$-skeleta of systolic complexes were earlier studied by metric graph theorists as \emph{bridged graphs} \cite{bandeltChepoi:metricGraphTheoryGeometrySurvey}.

    \begin{definition}
        Let $C\to X$ be an embedded cycle in a complex $X$. A $1$-cell $e$ with endpoints in $C$ is \emph{diagonal} if no edge of $C$ has the same endpoints as $e$. A simplicial complex is \emph{flag} if every finite set of pairwise adjacent $0$-simplices spans a simplex. A flag complex $\Delta$ is \emph{$6$-large} if every embedded cycle $C\to \Delta$ with $4\leq |C|<6$ has a diagonal. 
    \end{definition}

    \begin{definition}
        Let $\Delta$ be a simplicial complex, and let $x$ be a $0$-simplex of $\Delta$. The \emph{link} of $x$ is the simplicial complex $\lk(x)$ with a $k$-simplex for every corner of a $(k+1)$-simplex incident at $x$. A simply-connected, simplicial complex $\Delta$ is \emph{systolic} if $\lk(x)$ is $6$-large for every $0$-simplex $x$ of $\Delta$.
    \end{definition}

    \begin{definition}
        A \emph{tricomb} $\tflat$ is a $2$-complex isomorphic to the tiling of the Euclidean plane by regular triangles. A \emph{flat} in a systolic complex $\N$ is an isometric embedding $\tflat\to \N$. 
    \end{definition}

    \begin{definition}
        A systolic complex $\Delta$ has \emph{isolated flats} if for every $c>0$ there exists $K>0$ with $\diam\big(N_c(\tflat)\cap N_c(\mathbb{E'}^2_\mathbf{tri})\big)<K$ for distinct tricombs $\tflat$ and $\mathbb{E'}^2_\mathbf{tri}$ in $\Delta$.
    \end{definition}

    There is a rich literature on the theory of $\mathrm{CAT}(0)$ spaces with isolated flats. For an account of this literature, see the introduction of \cite{hruskaKleiner:hadamardSpacesWithIsolatedFlats}. In the same paper, Hruska and Kleiner prove groups acting properly and cocompactly on a $\mathrm{CAT}(0)$ space with isolated flats are relatively hyperbolic. Elsner proved the following systolic analogue of this theorem. 

    \begin{theorem}[\cite{elsner:systolicGroupsWithIsolatedFlats}]
    \label{thm:elsner}
        Suppose $G$ acts properly cocompactly on a systolic complex with isolated flats. Then $G$ is hyperbolic relative to a family of its maximal virtually abelian subgroups of rank $2$.
    \end{theorem}

    It suffices to specify the $0$-simplicies and $1$-simplicies when constructing a flag simplicial complex, since a flag simplicial complex is determined by its $1$-skeleton. 

    \begin{definition}
        Let $X$ be a $2$-complex. Its \emph{nerve} $\N X$ is the flag simplicial complex whose $0$-simplicies correspond to faces of $X$ and whose $1$-cells correspond to intersecting faces. A map $\phi:Y\to X$ \emph{respects isolation} if it sends isolated $1$-cells of $Y$ to isolated $1$-cells of $X$. For any $\phi:Y\to X$ that respects isolation, there is an induced map $\N \phi: \N Y\to \N X$.
    \end{definition}

    \begin{example}
        Two $2$-cells in a honeycomb intersect if and only if they share a hexside, so $\N \hflat = \tflat$.
    \end{example}

    The following follows from the definitions. 

    \begin{lemma}
    \label{lem:fdistNXisometric}
        Let $X$ be a $2$-complex. The natural map from the faces of $X$ with the $\fdist$-metric to $\N X$ is an isometry. We denote this map by $\N:X\to \N X$.
    \end{lemma}



    The following theorem was established by Wise \cite{wise:sixtolicComplexesAndTheirFundamentalGroups}. Later, Osajda and Prytu\l a generalized it to the graphical small-cancellation setting \cite[thm~7.12]{osajdaPrytula:classifyingSpacesForFamiliesofSubgroupsForSystolicGroups}.
    
    \begin{theorem}
    \label{them:wiseNerve}
        If $X$ is a simply-connected $\Csix$ complex, then $\N X$ is systolic. 
    \end{theorem}

    \begin{lemma}
    \label{lem:locallyFlat}
        Let $X$ be a simply-connected $\Csix$ complex. Let $x$ be a $0$-simplex in a flat $\tflat\subset \N X$. Let $y_1,\ldots, y_6$ be the neighbors of $x$ in $\tflat$. The subpaths $\N^{-1}y_i\cap \N^{-1}x$ of $\partial_c (\N^{-1}x)$ are nontrivial petal-pieces with associated petals $\N^{-1}y_i$.
    \end{lemma}

    \begin{proof}
        Let $R=\N^{-1}x$, $R_i=\N^{-1} y_i$, and $P_i=\N^{-1}y_i\cap \N^{-1} x=R_i\cap R$. By Corollary~\ref{cor:2cellIntersection}, each $P_i$ is a closed arc in $\partial R$, and $R_i$, $R_j$ intersect if and only if $P_i$, $P_j$ intersect. By Lemma~\ref{lem:fdistNXisometric}, $\N^{-1}:\N X\to X$ is an isometry. Hence for $i,j\in\Z/6\Z$, the intersection $P_i\cap P_j$ is nonempty if and only if $|i-j|=1$. Consequently, $P_1,\ldots, P_6$ are nontrivial pieces covering $\partial R$.      
    \end{proof}

    Before stating the next theorem, we emphasize a subtlety in the definitions of a flat. A flat in a combinatorial $2$-complex $X$ is an embedding $\hflat\to X$, which is not a priori isometrically embedded. However, a flat in a systolic complex $\N$ is an isometric embedding $\tflat\to \N$. 
    
    \begin{theorem}
    \label{thm:tricombToHoneycomb}
        Let $X$ be a simply-connected $\Csix$ complex. Then $\N:X\to \N X$ induces a bijection between flats in $X$ and flats in $\N X$.
    \end{theorem}

    \begin{proof}
        Let $F:\hflat\to X$ be a flat in $X$. By Lemma~\ref{lem:fdistNXisometric} and Lemma~\ref{lem:immersedHoneycombIsometric}, $\N F:\tflat \to \N X$ is a flat in $\N X$. It remains to show each flat $\tflat\to\N X$ is $\N F$ for some flat $F$ in $X$. 

        Identify the flat $\tflat\to \N X$ with its image $\tflat\subset \N X$, and consider the inverse image $\N^{-1}\tflat$ under $\N$. By Lemma~\ref{lem:locallyFlat}, $\N^{-1}\tflat$ is locally isomorphic to a honeycomb, so there exists a honeycomb $\hflat$ and an immersion $F:\hflat\to X$ which is a local isometry onto $\N^{-1}\tflat$. By Lemma~\ref{lem:immersedHoneycombIsometric}, $F$ is a flat in $X$, and $\N F$ is equal to the original flat $\tflat\subset \N X$. 
    \end{proof}

    Combining Lemma~\ref{thm:tricombToHoneycomb} with Lemma~\ref{lem:fdistNXisometric}, we obtain:

    \begin{corollary}
    \label{cor:isoTriIFFisoHex}
        Let $X$ be a simply-connected $\Csix$ complex. Then $X$ has isolated flats if and only if $\N X$ has isolated flats.
    \end{corollary}

    \begin{theorem}
    \label{thm:relativelyHyperbolic}
        If $G$ acts properly cocompactly on a simply-connected $\Csix$ complex $X$ with isolated flats, then $G$ is hyperbolic relative to a family of maximal virtually abelian subgroups of rank $2$. In particular, this holds if $X$ is strict $\Csix$.
    \end{theorem}

    \begin{proof}
    The map $\N:X\to \N X$ induces an action of $G$ on $\N X$. By Theorem~\ref{them:wiseNerve}, Corollary~\ref{cor:isoTriIFFisoHex}, and Theorem~\ref{thm:elsner}, $G$ is hyperbolic relative to a family of maximal virtually abelian subgroups of rank $2$. If $X$ is strict $\Csix$, then $X$ has isolated flats by Theorem~\ref{thm:strictC6isoFlats}.
    \end{proof}

\hfill
    
\section{Convex cocompact cores}
\label{sec:convexCores}

\subsection*{\S Cores in strict $\Csix$ complexes}

    \begin{definition}
        Let $A_1,\ldots, A_m, B$ and $B, C_1,\ldots, C_n$ be wall-segments such that $A_m\cap C_1=\emptyset$. The sequence $A_1,\ldots, A_m, B, C_1,\ldots, C_n$ is a \emph{bent wall-segment}.
    
    \end{definition}

    \begin{lemma}
    \label{lem:bentSegment}
        Let $X$ be a simply-connected strict $\Csix$ complex. Bent wall-segments in $X$ are $\fdist$-geodesics.
    \end{lemma}

    \begin{proof}
        Let $A_1,\ldots, A_m, B, C_1,\ldots, C_n$ be a bent wall-segment, and let $R_1,\ldots, R_k$ be a $\fdist$-geodesic joining $R_1=A_1$ and $R_k=C_n$. Let $J$ be a concatenation of two traces: one tracing $A_1,\ldots, B$, and another tracing $B,\ldots, C_n$. Suppose $J$ meets $A_1$ only at its initial vertex and $J$ meets $C_n$ only at its terminal vertex. Let $P$ be a trace of $R_1,\ldots, R_k$ with the same endpoints as $J$. By Lemma~\ref{lem:vanKampen}, there is a disc diagram $\phi:D\to X$ such that $\phi\circ (\partial_c D)=PJ^{-1}$. Choose $J$, $P$, and $D$ as above such that $D$ has a minimal number of cells. 

        Let $J=\alpha_2\cdots \alpha_m\beta\gamma_1\cdots \gamma_{n-1}$ where $\alpha_i\to A_i$ and $\beta\to B$ and $\gamma_i\to C_i$. We have $\plength{\alpha_i}>1$ and $\plength{\gamma_i}>1$ since each $\alpha_i$ and $\gamma_i$ joins pieces of $A_i$ and $C_i$ containing opposite $1$-cells. Consequently, no outerpath of a shell can contain an $\alpha_i$ or $\gamma_i$ as a subpath. Thus for any $2$-cell $R$ of $D$, we have $\plength{\partial_cR \cap J}\leq 3$. Moreover, if $\plength{\partial_cR\cap J}=3$ then $\partial_c R$ must contain $\beta$. In particular, $J$ contains no outerpath of a $2$-shell and at most one outerpath of a $3$-shell. By minimality, $\plength{\partial_cR\cap P}\leq 2$, otherwise $\partial_cR$ would meet three consecutive faces $R_{i-1}, R_i, R_{i+1}$ and could be removed from $D$, with $R$ replacing $R_i$ in the $\fdist$-geodesic $R_1,\ldots, R_k$. Thus $P$ contains no outerpath of a $3$-shell. By the above observations, $D$ contains at most three $3$-shells: one at each of the shared endpoints of $P$ and $J$, and one whose outerpath contains $\beta$. However, by Lemma~\ref{lem:greendlinger}, if $D$ contains exactly three shells, then $D$ contains three $2$-shells, a contradiction, since $J$ cannot contain the outerpath of a $2$-shell. Hence $D$ must be a ladder $L_1,\ldots, L_k$, which we claim is a ladder of $1$-cells so that $J=P$. Indeed, suppose $L_i$ is the first $2$-cell of $D$.As observed above, $\plength{\partial_c L_i\cap J}\leq 3$ and $\plength{\partial_c L_i\cap P}\leq 2$. By $\Csix$, we must have $\plength{\partial_c L_i\cap J}=3$ and $\plength{\partial_c L_i\cap P}\leq 2$ and $L_{i+1}$ is a $2$-cell forming a piece with $L_i$. Furthermore, $L_i$ contains $\beta$. Thus $\plength{\partial_c L_{i+1}\cap J}<3$ and $\plength{\partial_c L_{i+1}\cap P}=2$, so there must exist a $2$-cell $L_{i+2}$ in the ladder forming a piece with $L_{i+1}$. This sequence never terminates, a contradiction.
        
        Using $J=P$, we consider how the decompositions $\alpha_2\cdots \alpha_m\beta\gamma_1\cdots \gamma_{n-1}$ and $P_1\cdots P_k$ overlap. No $P_i$ contains $\alpha_i$ since $\plength{\alpha_i}>1$ for $i\in\{2,\ldots, m\}$. Similarly, no $P_{i+m}$ contains $\gamma_i$ since $\plength{\gamma_i}>1$ for $i\in\{1,\ldots, n\}$. Consequently, $k\geq m+n-1$ and thus the bent wall-segment $A_1,\ldots, A_m, B, C_1,\ldots, C_n$ is a $\fdist$-geodesic. 
    \end{proof}

    The following definition will be used in the proof of Theorem~\ref{thm:convexHull}  

\begin{figure}[h]
        \centering
        \includegraphics[width=0.6\textwidth]{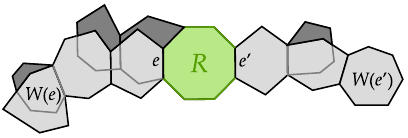}
        \caption{The green indicates the cells to be deleted, leaving the components $W(e)$ and $W(e')$.}	
        \label{fig:wallray}
    \end{figure}

    \begin{definition}
        Let $W(e,e')$ be a carrier containing opposite $1$-cells $e$ and $e'$. Remove each $2$-cell $R$ of $W(e,e')$ containing $e$, $e'$ and each $1$-cell of $R$ not contained in a piece with either $e$ or $e'$. See Figure~\ref{fig:wallray}. This separates $W(e,e')$ into two components. Let $W(e)$ denote the component containing $e$, and let $W(e')$ denote the component containing $e'$. 
    \end{definition}

    \begin{definition}
        Let $X$ be a simply-connected strict $\Csix$ complex. A subcomplex $Y\subset X$ is \emph{$k$-quasiconvex} if $R_1,\ldots, R_n\subset N_k(Y)$ whenever $R_1,\ldots, R_n$ is a $\fdist$-geodesic and $R_1\cap Y\neq\emptyset$, $R_n\cap Y\neq \emptyset$. We say $Y$ is \emph{quasiconvex} if it is $k$-quasiconvex for some $k$.
    \end{definition}

    \begin{definition}[$\Hull(Y)$]
        Let $Y$ be a subcomplex of a simply-connected strict $\Csix$ complex. The \emph{hull of $Y$} is the minimal face-convex subcomplex containing $Y$. 
    \end{definition}
    
    \begin{theorem}
    \label{thm:convexHull}
        Let $X$ be a simply-connected strict $\Csix$ complex, and let $Y$ be a quasiconvex subcomplex. For some $K>0$, we have $\Hull(Y)\subset N_K(Y)$. 
    \end{theorem}

    \begin{proof}
        Let $k$ be a constant such that any $\fdist$-geodesic $R_1,\ldots, R_n$ with $R_1\cap Y\neq \emptyset$ and $R_n\cap Y\neq \emptyset$ is contained in $N_k(Y)$. Let $K>k+2$. We claim if $R$ is a face with $\fdist(R,Y)>K$, then there exists a wall $\wall$ separating $R$ and $Y$. That is, $R$ and $Y$ lie in distinct halfspaces of $\wall$. The theorem follows from this claim since $\Hull(Y)$ is contained in the intersection of halfspaces containing $Y$.

        Let $R$ be a face with $\fdist(R,Y)>K$. Let $R'_1,\ldots, R'_n=R$ be a $\fdist$-geodesic with $R'_1\cap Y\neq \emptyset$. 
        Let $e', f'$ be opposite $1$-cells of $R'_{n-1}$ such that any piece $P\to R'_{n-1}$ containing $e'$ or $f'$ is disjoint from $R'_n$. Let $W(e',f')$ be a carrier of a wall containing $e'$, $f'$. The carrier $W(e',f')$ contains $R'_{n-1}$ but not $R'_n$ by our choice of $e'$ and $f'$. If $W(e',f')$ is disjoint from $Y$, then it separates $R$ and $Y$. Otherwise, $W(e',f')$ intersects $Y$. Let $R_1,\ldots, R_m=R'_{n-1}$ be a wall segment in $W(e',f')$ with $R_1\cap Y\neq \emptyset$. 

        Let $e, f$ be the $1$-cells of $R_{m-1}$ closest to $R_{m-1}\cap R_{m-2}$ such that any piece containing $e$ or $f$ is disjoint from $R_{m-1}\cap R_{m-2}$. There are now two cases depending on whether $e$, $f$ are opposite.
        
        Suppose that $e$ and $f$ are opposite, and consider a carrier $W(e,f)$. We claim that $W(e,f)$ separates $R$ and $Y$. Indeed, $R_{m-1}, R_m, R$ is a $\fdist$-geodesic by Lemma~\ref{lem:bentSegment}, and $R_m\not\subset W(e,f)$. Thus $R\not\subset W(e,f)$ by Lemma~\ref{lem:wallSegmentsInWalls}. Also $W(e,f)$ does not intersect $Y$, otherwise there would exist a wall-segment $R_{m-1},B_l,\ldots, B_1$ with $B_1\cap Y\neq\emptyset$. But then $R_1,\ldots, R_{m-1},\ldots, B_1$ is a bent wall-segment, and thus a $\fdist$-geodesic by Lemma~\ref{lem:bentSegment}. Since $R_1,B_1$ intersect $Y$, $\fdist(R_{m-1},Y)<k$, contradicting our choice of $K$. Hence $W(e,f)$ is disjoint from both $R$ and $Y$, so it separates $R$ and $Y$. 

        Suppose $e, f$ are not opposite. Then there exists an embedded path which is the concatenation of three pieces $P_1P_2P_3\to \partial R_{m-1}$ such that $e$ is the initial edge of $P_1$ and $f$ is the terminal edge of $P_3$. By the choice of $e$ and $f$, the path $P_1P_2P_3$ is disjoint from $R_{m-1}\cap R_{m-2}$. Moreover, there exists a petal piece decomposition $\partial_cR_{m-1}=S_1S_2S_3P_1P_2P_3$ where $S_2$ parametrizes $R_{m-1}\cap R_{m-2}$. Let $e'$ be the initial $1$-cell of $S_1$, and let $f'$ be the terminal $1$-cell of $S_3$. By Lemma~\ref{lem:petalPieceMaximal}, $e,e'$ and $f,f'$ are pairs of opposite $1$-cells in $R_{m-1}$, and we let $W(e,e')$ and $W(f,f')$ be carriers containing these pairs. 
        Neither $W(e)$ nor $W(f)$ can intersect $Y$. Indeed, if $R_{m-1},B_l,\ldots, B_1$ is a wall-segment in either $W(e)$ or $W(f)$ with $B_1\cap Y\neq \emptyset$, then $R_1,\ldots, R_{m-1}, B_l, \ldots, B_1$ is a bent wall-segment whose ending faces intersect $Y$, so $\fdist(R_{m-1},Y)<k$, a contradiction. A similar argument shows that at most one of $W(e')$ or $W(f')$ can intersect $Y$. Without loss of generality, suppose $W(e')$ does not intersect $Y$. Observe $R\not\subset W(e,e')$ since $R_{m-1},R_m,R$ is a $\fdist$-geodesic by Lemma~\ref{lem:bentSegment} with $R_m\not\subset W(e,e')$. Since $W(e,e')$ is also disjoint from $Y$, we conclude that $W(e,e')$ separates $Y$ and $R$.
    \end{proof}

    \begin{corollary}
    \label{cor:convexCore}
        Let $G$ be a hyperbolic group acting properly cocompactly on a simply-connected strict $\Csix$ complex $X$. Let $H$ be a quasiconvex subgroup of $G$, and let $C$ be a compact subcomplex of $X$. There exists a $\fdist$-convex $H$-cocompact subcomplex $Y$ of $X$ containing $HC$. 
    \end{corollary}

    Corollary~\ref{cor:convexCore} explains that quasiconvexity in a hyperbolic strict $\Csix$ group is witnessed by $\fdist$-convexity. We describe the relatively quasiconvex analogue in Theorem~\ref{thm:cosparseCore}. 

    \begin{remark}
    \label{rem:relativelyQuasiconvex}
        Suppose $G$ acts properly cocompactly on a simply-connected strict $\Csix$ complex $X$. By Theorem~\ref{thm:relativelyHyperbolic}, $G$ is hyperbolic relative to a family of maximal virtually abelian subgroups of rank $2$. There are a number of equivalent definitions of a relatively quasiconvex subgroup \cite{hruska:relativeHyperbolicityRelativeQuasiconvexity}. The definition we give below works in this context by \cite[cor~8.16]{hruska:relativeHyperbolicityRelativeQuasiconvexity}. 
    \end{remark}
    
    \begin{definition}
    \label{def:relativelyQuasiconvex}
        Let $G$ and $X$ be as in Remark~\ref{rem:relativelyQuasiconvex}. Let $H$ be a subgroup of $G$, and let $\mathcal F$ be the set of honeycombs $\hflat$ such that $\Stab(\hflat)\cap H$ is infinite. We say $H$ is \emph{relatively quasiconvex} if $HC\cup \bigcup_{F\in \mathcal F} F$ is a quasiconvex subcomplex for some (hence any) compact subcomplex $C$ of $X$.
    \end{definition}

    The following is a consequence of \cite[thm~10.5]{hruska:relativeHyperbolicityRelativeQuasiconvexity}.

    \begin{lemma}
    \label{lem:quasiIsometricEmbedding}
        Let $G$ and $X$ be as in Remark~\ref{rem:relativelyQuasiconvex}. Let $H$ be a relatively quasiconvex subgroup of $G$, and let $R$ be a face of $X$. The orbit map $H\to X:h\mapsto hR$ is a quasi-isometric embedding.
    \end{lemma}

    We will make use of relatively thin triangles. The following theorem is a restatement of \cite[thm~4.1]{sageevWise:coresForQuasiconvexActions}, which is extracted from \cite[cor~8.14, cor~8.19]{drutuSapir:treeGradedSpacesAndAsymptoticCones}.

\begin{figure}[h]
        \centering
        \includegraphics[width=0.5\textwidth]{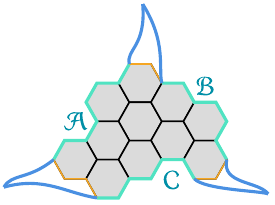}
        \caption{The subgments $A$, $B$, $C$ of $\mathcal A$, $\mathcal B$, $\mathcal C$ contained in $\hflat$ are highlighted. The orange segments join endpoints of $A$, $B$ and $B$, $C$ and $C$, $A$. Their lengths are uniformly bounded by $K$.}	
        \label{fig:triangle}
    \end{figure}

    \begin{theorem}
    \label{thm:relativelyThinTriangles}
    Let $G$ and $X$ be as in Remark~\ref{rem:relativelyQuasiconvex}. For each $k>0$ there exists $K>0$ so that for any $\fdist$ $k$-quasigeodesic triangle $\mathcal A\mathcal B\mathcal C$ in $X$ one of the following holds:
        \begin{enumerate}
            \item There exists a face $R$ such that $N_{K/2}(R)$ intersects each of $\mathcal A$, $\mathcal B$, $\mathcal C$. 
            \item There exists a honeycomb $\hflat$ such that $N_K(\hflat)$ contains subsegments $A$, $B$, $C$ of $\mathcal A$, $\mathcal B$, $\mathcal C$. Furthermore, the endpoints of $A$, $B$ closest to the corner formed by $\mathcal A$, $\mathcal B$ are within $K$ of each other. Similar statements hold for $B$,~$C$ and $C$, $A$. See Figure~\ref{fig:triangle}.
        \end{enumerate}
    \end{theorem}

    \begin{definition}
        A \emph{thick flat} $\mathbf E$ is a subcomplex containing a honeycomb $\hflat$ such that $\mathbf E$ and $\hflat$ are at a finite Hausdorff distance.
    \end{definition}

    \begin{definition}
    \label{def:cosparse}
        A proper action of $H$ on a simply-connected $\Csix$ complex $Y$ is \emph{cosparse} if there is a compact subcomplex $C$ and a finite collection of thick flats $\{\mathbf F_i\}$ such that the following hold:
        \begin{itemize}
            \item $Y=HC\cup \bigcup_i\mathbf F_i$.
            \item $Y\subset N_K(\Hull(HC))$ for some $K>0$.
            \item For each $i,j$ and $h\in H$, either $(\mathbf F_i\cap h\mathbf F_j)\subset HC$ or $\mathbf F_i=h\mathbf F_j$.
        \end{itemize}
    \end{definition}

    \begin{theorem}
    \label{thm:cosparseCore}
        Suppose $G$ acts properly cocompactly on a simply-connected $\Csix$ complex $X$. Let $H$ be a relatively quasiconvex subgroup of $G$, and let $C$ be a compact subcomplex of $X$. There exists a $\fdist$-convex $H$-cosparse subcomplex $Y$ containing $HC$.
    \end{theorem}

    \begin{proof}
        We may assume that $C$ contains a face $R$ of $X$. Let $\mathcal F$ be as in Definition~\ref{def:relativelyQuasiconvex}, and set $Y=\Hull(HC\cup \bigcup_{F\in\mathcal F}F)$. Since $H$ is relatively quasiconvex, $HC\cup \bigcup_{F\in\mathcal F}F$ is a quasiconvex subcomplex. Since $H$ is relatively quasiconvex, $\mathcal F$ is a union of finitely-many $H$-orbits of honeycombs. Let $\{F_i\}$ be a set of orbit representatives. By Theorem~\ref{thm:convexHull}, $Y\subset N_K(HC\cup \bigcup_iHF_i)$ for some $K>0$. The subcomplex $\mathbf C=Y\cap N_K(C)$ is compact, and each $\mathbf F_i=Y\cap N_K(F_i)$ is a thick flat. We have $Y=H\mathbf C\cup\bigcup_iH\mathbf F_i$, so $Y$ satisfies the first condition of Definition~\ref{def:cosparse}. We now adjust $Y$ so that it satisfies the second condition.

        Enlarge $\mathbf C$ by a finite amount so that $H\mathbf C\cap F_i$ is nonempty for each $i$. The intersection $L=H\mathbf C\cap F_i$ is either coarsely dense in $F_i$ or a quasiline. Suppose $L$ is finite Hausdorff distance from a width $1$ band $B\subset F_i$. A wall carrier in $\mathbf F_i$ containing a band $B'\subset F_i$ sufficiently far and parallel to $B$ is disjoint from $L$.  Hence $\Hull(L)\cap \mathbf F_i$ is coarsely equal to $L$. If $L$ is not finite Hausdorff distance from a width $1$ band in $F_i$, then consideration of bent wall-segments shows $F_i\subset \Hull(L)$. See Figure~\ref{fig:diagonalQuasiline}. Remove each thick flat in $\{\mathbf F_i\}$ for which $\Hull(L)\cap \mathbf F_i$ is coarsely equal to $L$, and enlarge $\mathbf C$ by a finite amount so that $H\mathbf C$ contains $\Hull(L)\cap \mathbf F_i$ for each removed $\mathbf F_i$. Replace each remaining thick flat $\mathbf F_j$ with $\Hull(L)\cap \mathbf F_j$. Then we have $\Hull(L)\cap \mathbf F_j\subset \Hull(H\mathbf C)$, so that $Y$ satisfies the second condition of Definition~\ref{def:cosparse}.

        \begin{figure}[h]
            \centering
            \includegraphics[width=0.5\textwidth]{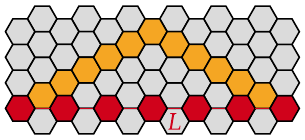}
            \caption{The red subcomplex $L=H\mathbf C\cap F_i$ is a quasiline not parallel to any width $1$ band. The orange subcomplex intersecting two distant faces of $L$ is a bent wall-segment and thus a geodesic by Lemma~\ref{lem:bentSegment}.}	
        \label{fig:diagonalQuasiline}
        \end{figure}

        We claim that $\mathbf F_i\cap h\mathbf F_j$ is contained in a uniformly finite neighborhood of $H\mathbf C$ for each $i$, $j$ and $h\in H$. From this claim, after potentially enlarging $\mathbf C$ by a finite amount, it follows that $Y$ satisfies the third condition of Definition~\ref{def:cosparse}. Let $R$ be a face of $\mathbf F_i$ intersecting $h\mathbf F_j$. Since $H\mathbf C\cap \mathbf F_i$, $H\mathbf C\cap h\mathbf F_j$ are nonempty and $F_i$, $hF_j$ are isometrically embedded, there exist $\fdist$-quasigeodesics in $\mathcal A$, $\mathcal B$ in $\mathbf F_j$, $h\mathbf F_j$ joining $R$ to $H\mathbf C$. By Lemma~\ref{lem:quasiIsometricEmbedding}, the ends of $\mathcal A$, $\mathcal B$ in $H\mathbf C$ are joined by a quasigeodesic $\mathcal C$ which stays in a uniformly finite neighborhood of $H\mathbf C$. We apply Theorem~\ref{thm:relativelyThinTriangles} to $\mathcal A\mathcal B\mathcal C$ and consider the two possible cases. If $N_{K/2}(R')$ intersects each of $\mathcal A$, $\mathcal B$, $\mathcal C$ as in Theorem~\ref{thm:relativelyThinTriangles}(1), then $R$ is uniformly close to $H\mathbf C$ since coarse intersections of honeycombs are bounded by Theorem~\ref{thm:strictC6isoFlats}. Suppose $N_K(\hflat)$ contains subsegments $A$, $B$, $C$ of $\mathcal A$, $\mathcal B$, $\mathcal C$ as in Theorem~\ref{thm:relativelyThinTriangles}(2). Let $A'$ be the subsegment of $\mathcal A\setminus A$ containing $R$. The endpoints of $A'$ lie in the coarse intersection of two honeycombs, and $A$ is contained in the coarse intersection of honeycombs. Thus both $A$ and $A'$ have uniformly bounded length by Theorem~\ref{thm:strictC6isoFlats}. The path $A'A$ joins $R$ to a uniformly finite neighborhood of $H\mathbf C$ and has uniformly bounded length, so this concludes the proof.
    \end{proof}

\subsection*{\S Cores in non-strict $\Csix$ complexes}

Our goal here is to describe failures of convex cocompact cores. We first give an example of a $\Csix$ complex using the graph metric. We then describe a more interesting $\mathrm{C}(4)$ example using the $\fdist$-metric. 

\begin{example}
    The convex hull of $\{0,1,2,3\}$ is the entire Cayley graph of $\Z=\langle 1,4\rangle$. Indeed, $4$ lies on a geodesic from $0$ to $3$, and this repeats in both directions. We can make this a $\Csix$ example by using the presentation $\Z=\langle a,d\mid d=a^4\rangle$. The $2$-cell in this presentation has a free face so the complex is strict $\Csix$. Hence this complex has the convex cocompact core property with respect to the $\fdist$-metric but not the graph metric.
\end{example}

\begin{figure}[h]
    \centering
    \includegraphics[width=0.2\textwidth]{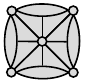}
    \caption{The $2$-complex $\mathbf S$ has seventeen $2$-cells. Eight pentagons surround the central octagon. The subcomplex $\mathbf B$ consists of the eight boundary triangles.}	
    \label{fig:thickSquare}
\end{figure}

\begin{definition}
    Let $\mathbf S$ be the $2$-complex depicted in Figure~\ref{fig:thickSquare}. The subcomplex $\mathbf B\subset \mathbf S$ is the union of the boundary $2$-cells. 
\end{definition}

\begin{remark}\label{rem:isometricBoundary}
    Note $\mathbf B\subset\mathbf S$ is $\fdist$-isometrically embedded, and $\Hull(\mathbf B)=\mathbf S$.
\end{remark}

\begin{figure}[h]
    \centering
    \includegraphics[width=0.55\textwidth]{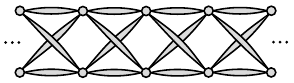}
    \caption{The complex $\mathbf X'$ obtained by blowing-up a Cayley graph of $\Z\times \Z_2$.}	
    \label{fig:thickCayley}
\end{figure}

\begin{definition}
    Consider the Cayley graph of $\Z\times \Z_2$ with generating set $\{(1,0),(1,1)\}$. Blow-up its $0$-cells and $1$-cells to discs to obtain the $2$-complex $\mathbf X'$ depicted in Figure~\ref{fig:thickCayley}.
\end{definition}

\begin{example}\label{ex:unboundedHull}

    \begin{figure}[h]
        \centering
        \includegraphics[width=0.75\textwidth]{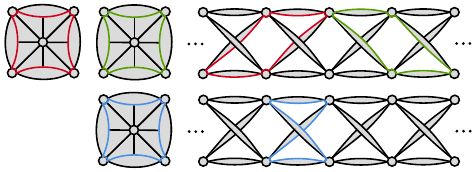}
        \caption{The blown-up $4$-cycles in $\mathbf X'$ containing red and green edges are two translates of $\mathbf C_\alpha$. The blown-up $4$-cycle containing blue edges is a translate of $\mathbf C_\beta$. The edge colorings indicate how copies of $\mathbf S$ are attached.}	
        \label{fig:fourCycle}
    \end{figure}
    
    The action of $\mathbf Z\times \mathbf Z_2$ on its Cayley graph induces an action on $\mathbf X'$. The vertices $\{(0,0), (1,0), (2,1), (1,1)\}$ span a $4$-cycle in the Cayley graph, and their blow-up is a subcomplex $\mathbf C_\alpha\subset \mathbf X'$. Up to subdividing $1$-cells, the subcomplexes $\mathbf C_\alpha$ and $\mathbf B$ are isomorphic. Equivariantly attach a copy of $\mathbf S$ along $\mathbf B$ to each translate of $\mathbf C_\alpha$, subdividing as necessary. See Figure~\ref{fig:fourCycle}. The $4$-cycle $\{(0,0),(1,0),(0,1),(1,1)\}$ in the Cayley graph blows-up to a subcomplex $\mathbf C_\beta$. Equivariantly attach copies of $\mathbf S$ along $\mathbf B$ to translates of $\mathbf C_\beta$, subdividing $1$-cells as needed. The action of $\Z\times \Z_2$ extends to the resulting complex $\mathbf X$, which is simply-connected. Note that $\Z_2$ exchanges the red and green cycles while preserving the blue cycles in Figure~\ref{fig:fourCycle}.

    Let $R\subset \mathbf X$ be the $2$-cell blow-up of a vertex in the Cayley graph of $\Z\times \Z_2$. Let $R'=(0,1)R$ so that $\Z_2R=R\cup R'$. We claim that the convex hull of $\mathbf Z_2R$ is all of $\mathbf X$. One can verify that $\fdist(R,R')=4$. Thus $\Hull(\mathbf Z_2R)$ contains $\mathbf B\subset \mathbf S$ for each copy of $\mathbf S$ containing $R$ and $R'$. See Figure~\ref{fig:hullStep1}. In particular, $\Hull(\mathbf Z_2R)$ contains the $(1,0)$ and $(-1,0)$-translates of $\mathbf Z_2R$. Repeating the above arguments for these translates of $\mathbf Z_2R$, we get that $\mathbf X'\subset \Hull(\mathbf Z_2R)$. By Remark~\ref{rem:isometricBoundary}, $\Hull(\mathbf X')=\mathbf X$. 

    \begin{figure}[h]
        \centering
        \includegraphics[width=0.55\textwidth]{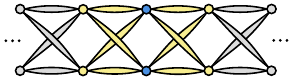}
        \caption{The $2$-cells $R$ and $R'$ are blue. The yellow (and blue) $2$-cells are the copies of $\mathbf B$ containing $R$ and $R'$. Note the yellow cells are the union of four $\fdist$-geodesics joining $R$ and $R'$ and thus contained in $\Hull(\mathbf Z_2R)$.}	
        \label{fig:hullStep1}
    \end{figure}
\end{example}

\bibliographystyle{alpha}
\bibliography{bibliography.bib}
    
\end{document}